\documentclass[11pt]{article}
\usepackage{epsfig}
\usepackage{amssymb}
\usepackage{amscd}
\usepackage{epsf}
\usepackage{amsfonts,amssymb,amsthm,amsmath}
 \usepackage{graphicx}
\usepackage[matrix,arrow,curve]{xy}

\input cyracc.def

\long\def\comment#1\endcomment{\relax}

\newcounter{subsubsubsection}
\newcounter{subsubsubsubsection}

\makeatletter \@addtoreset{subsubsubsection}{subsubsection}
\@addtoreset{subsubsubsubsection}{subsubsubsection}

\makeatother

\newcommand{\sevafigc}[4]{\begin{figure}[h]\centerline{
 \epsfig{file=#1,width=#2,angle=#3}}
\bigskip\caption{#4}\end{figure}}

\DeclareMathSizes{11.1}{10}{8}{6}

\DeclareMathOperator{\Hom}{Hom}

\newtheorem*{theorem*}{Theorem}

\newtheorem*{lemma}{Lemma}

\newtheorem*{klemma}{Key-Lemma}

\newtheorem*{proposition}{Proposition}

\newtheorem*{corollary}{Corollary}

\theoremstyle{remark}

\newtheorem*{remark}{Remark}

\newtheorem*{example}{Example}

\theoremstyle{definition}

\newtheorem*{definition}{Definition}

\newcommand{\hdot}{{\:\protect\raisebox{3pt}{\text{\protect\circle*{1.5}}}}}

\newcommand{\mb}{\hdot}

\newcommand{\RHom}{\mathrm{RHom}}

\newcommand\Lie{\mathrm{Lie}}

\newcommand\cat{\mathsf{cat}}
\newcommand\Free{\mathrm{Free}}

\newcommand{\Mor}{\mathrm{Mor}}
\newcommand{\Ext}{\mathrm{Ext}}

\newcommand{\gl}{\mathfrak{gl}}

\newcommand{\U}{\mathcal{U}}

\newcommand{\Assoc}{\mathrm{Assoc}}

\newcommand{\g}{\mathfrak{g}}

\newcommand{\Mod}{\mathrm{Mod}}
\newcommand{\Comm}{\mathrm{Comm}}
\newcommand{\Aff}{\mathrm{Aff}}
\newcommand{\gr}{\mathrm{gr}}

\newcommand{\Koszul}{\mathrm{Koszul}}

\newcommand{\Hoch}{\mathrm{Hoch}}
\newcommand{\poly}{{poly}}

\newcommand{\Log}{\mathrm{Log}}

\newcommand{\tot}{\mathrm{tot}}
\newcommand{\F}{\mathcal{F}}

\newcommand{\Ob}{\mathrm{Ob}}
\newcommand{\opp}{\mathrm{opp}}
\newcommand{\Rhom}{\mathrm{RHom}}
\newcommand{\Cone}{\mathrm{Cone}}

\newcommand{\Pois}{\mathrm{Pois}}

\def\wtilde#1{\widetilde{#1}\vphantom{#1}}

\title{ {\huge  Koszul duality in deformation quantization and Tamarkin's approach to Kontsevich formality}}
\author{{\LARGE Boris Shoikhet}}
\date{}

\begin{document}\maketitle

\begin{abstract}
Let $\alpha$ be a quadratic Poisson bivector on a vector space $V$.
Then one can also consider $\alpha$ as a quadratic Poisson bivector
on the vector space $V^*[1]$. Fixed a universal deformation
quantization (prediction of some complex weights to all Kontsevich graphs
[K97]), we have deformation quantization of the both algebras
$S(V^*)$ and $\Lambda(V)$. These are graded quadratic algebras, and
therefore Koszul algebras. We prove that for some universal
deformation quantization, independent on $\alpha$,  these two algebras are Koszul dual. We
characterize some deformation quantizations for which this theorem
is true in the framework of the Tamarkin's theory [T1].
\end{abstract}

\section*{Introduction}
This paper is devoted to the theorem that there exists a universal
deformation quantization compatible with the Koszul duality, as it
is explained in the Abstract. Let us firstly formulate it here in a
more detail, and then outline the main ideas of the proof.

\subsection{}
Let $V$ be a finite-dimensional vector space over
$\mathbb{C}$. We denote by $T_\poly(V)$ the graded Lie algebra of
polynomial polyvector fields on $V$, with the Schouten-Nijenhuis
bracket. For a $\mathbb{Z}$-graded vector space $W$ denote by $W[1]$
the graded space for which $(W[1])^i=W^{1+i}$, that is, the space
"shifted to the left". The following simple statement is very
fundamental for this work:
\begin{lemma}
There is a canonical isomorphism of graded Lie algebras
$\mathcal{D}\colon T_\poly(V)\to T_\poly(V^*[1])$.
\end{lemma}
The map $\mathcal{D}$ maps a bi-homogeneous polyvector field
$\gamma$ on $V$, $\gamma=x_{i_1}\dots
x_{i_k}\frac{\partial}{\partial
x_{j_1}}\wedge\dots\wedge\frac{\partial}{\partial x_{j_\ell}}$ to
the polyvector field
$\mathcal{D}(\gamma)=\xi_{j_1}\dots\xi_{j_\ell}\frac{\partial}{\partial
\xi_{i_1}}\wedge\dots\wedge\frac{\partial}{\partial \xi_{i_k}}$ on
the space $V^*[1]$. Here $\{x_i\}$ is a basis in $V^*$, and
$\{\xi_i\}$ is the dual basis in $V[-1]$.\qed

It is a good place to recall the Hochschild-Kostant-Rosenberg
theorem by which the cohomological Hochschild complex of the algebra
$A=S(V^*)$ endowed with the Gerstenhaber bracket has cohomology
isomorphic to $T_\poly(V)$ as graded Lie algebra. That is, the
Gerstenhaber bracket becomes the Schouten-Nijenhuis bracket on the
level of cohomology.

This theorem is related to the lemma above (which is certainly clear
just straightforwardly, without this more sophisticated argument),
as follows. Consider the algebras $A=S(V^*)$ and
$B=\Lambda(V)=Fun(V^*[1])$. The algebras $A$ and $B$ are Koszul dual
(see, e.g. [PP]). Bernhard Keller proved in [Kel1] (see also the
discussion below in Sections 1.4-1.6) that the cohomological
Hochschild complexes $\Hoch^\mb(A)$ and $\Hoch^\mb(B)$ are
quasi-isomorphic with all structures when $A$ and $B$ are quadratic
Koszul and Koszul dual algebras. (For the Hochschild cohomology it
was known before, see the references in loc.cit). In our case
$HH^\mb(A)=T_\poly(V)$ and $HH^\mb(B)=T_\poly(V^*[1])$.

\subsection{} The isomorphism $\mathcal{D}$ from the lemma above
does not change the grading of the polyvector field, but it maps
$i$-polyvector fields with $k$-linear coefficients to $k$-polyvector
fields with $i$-linear coefficients. In particular, it maps
quadratic bivector fields on $V$ to quadratic bivector fields on
$V^*[1]$. Moreover, $\mathcal{D}$ maps a Poisson quadratic bivector
on $V$ to a Poisson quadratic bivector on $V^*[1]$, because it is a
map of Lie algebras.

In [K97], Maxim Kontsevich gave a formula for deformation
quantization of algebra $S(V^*)$ by a Poisson bivector $\alpha$ on
$V$ (the vector spaced $V$ should be finite-dimensional). His
formula is organized as a sum over admissible graphs, and each graph
is taken with the Kontsevich weight $W_\Gamma$. In particular, this
$W_\Gamma$ depends only on the graph $\Gamma$ and does not depend on
dimension of the space $V$.

Consider now all these complex numbers $W_\Gamma$ as undefined
variables. Then the associativity gives an infinite number of
quadratic equations on $W_\Gamma$. Kontsevich's paper [K97] then
shows that these equations have at least one complex solution.
Actually there is a lot of essentially different solutions, as is
clear from [T] (see the discussion in Section 3 of this paper). Any
such deformation quantization is called {\it universal} because the
complex numbers $W_\Gamma$ do not depend on the vector space $V$.

The case of a quadratic Poisson bivector $\alpha$ is distinguished,
by the following lemma:
\begin{lemma}
Let $S(V^*)_\alpha$ be a universal deformation quantization of
$S(V^*)$ by a quadratic Poisson bivector $\alpha$. Then the algebra
$S(V^*)_\alpha$ is graded. This means that for $f\in
S^i(V^*)[[\hbar]]$ and $g\in S^j(V^*)[[\hbar]]$, the product $f\star
g\in S^{i+j}(V^*)[[\hbar]]$.
\begin{proof}
By the Kontsevich formula,
\begin{equation}\label{eq10.1}
f\star g=f\cdot g+\sum_{k\ge 1}\hbar^k\sum_{\Gamma\in
G_{k,2}}W_\Gamma B_\Gamma(f,g)
\end{equation}
where $G_{k,2}$ is the set of admissible graphs with two vertices on
the real line and $k$ vertices in the upper half-plane (see [K97],
Section 1 for details). Now each graph $\Gamma$ from $G_{k,2}$ has
$k$ vertices at the half-plane, and $2k$ edges. One can compute the
grading degree of $B_\Gamma(f,g)$ as follows. It is the sum of
degrees of quantities associated with all vertices (which is $\deg
f+\deg g +k\deg \alpha=\deg f+\deg g +2k$) minus the number of edges
(equal to $2k$ by definition of an admissible graph) because each
edge differentiate once, and then decreases the degree by 1). This
difference is equal to $\deg f+\deg g$.
\end{proof}
\end{lemma}
In particular, for quadratic deformation quantization the map
$x_i\cdot x_j\mapsto x_i\star x_j$ gives a
$\mathbb{C}[[\hbar]]$-linear endomorphism of the space
$S^2(V)[[\hbar]]$ which is clearly non-degenerate. We can find an
inverse to it, then we can present the star-algebra as the quotient
of the tensor algebra $T(V^*)$ by the set of {\it quadratic}
relations $R_{ij}\in V^*\otimes V^*$, one relation for each pair of
indices $1\le i<j\le \dim V$. We conclude, that the Kontsevich
deformation quantization of $S(V^*)$ by a quadratic Poisson bivector
is a graded quadratic algebra.

\subsection{} We actually get two
quadratic associative algebras for any universal deformation
quantization, one is the deformation quantization of $S(V^*)$ by the
quadratic Poisson bivector $\alpha$, and another one is the
deformation quantization of $\Lambda(V)=Fun(V^*[1])$ by the
quadratic Poisson bivector $\mathcal{D}(\alpha)$. Denote these two
algebras by $S(V^*)\otimes\mathbb{C}[[\hbar]]_\alpha$ and
$\Lambda(V)\otimes\mathbb{C}[[\hbar]]_{\mathcal{D}(\alpha)}$.

In the present paper we prove the following result:
\begin{theorem*}
There exists a universal deformation quantization such that the two
algebras $S(V^*)\otimes\mathbb{C}[[\hbar]]_\alpha$ and
$\Lambda(V)\otimes\mathbb{C}[[\hbar]]_{\mathcal{D}(\alpha)}$ are
Koszul dual as algebras over $\mathbb{C}[[\hbar]]$. In particular,
\begin{equation}\label{eqnewintro1}
\Ext^\mb_{S(V^*)\otimes\mathbb{C}[[\hbar]]_\alpha-Mod}(\mathbb{C}[[\hbar]],\mathbb{C}[[\hbar]])=\Lambda(V)\otimes\mathbb{C}[[\hbar]]_{\mathcal{D}(\alpha)}
\end{equation}
and
\begin{equation}\label{eqnewintro2}
\Ext^\mb_{\Lambda(V)\otimes\mathbb{C}[[\hbar]]_{\mathcal{D}(\alpha)}-Mod}(\mathbb{C}[[\hbar]],\mathbb{C}[[\hbar]])=S(V^*)\otimes\mathbb{C}[[\hbar]]_\alpha
\end{equation}
The Tamarkin's deformation quantization defined from any Drinfeld's
associator (which is clearly universal) satisfies the condition of
Theorem.
\end{theorem*}
See Section 1 of this paper for an overview of Koszul duality, and
of Koszul duality over a discrete valuation ring.

\begin{remark}
To consider $S(V^*)$ and $\Lambda(V)$ as Koszul dual algebras, the $\Ext$ groups above should be taken in the category of $\mathbb{Z}$-graded modules over $\mathbb{Z}$-graded algebras. Without considering the $\mathbb{Z}$-graded category, the Koszul dual to $\Lambda(V)$ is $S[[V]]$. It is everywhere implicitly assumed that we work in the $\mathbb{Z}$-graded category.
\end{remark}

\subsection{}\label{section0.4}
Now let us outline our strategy how to prove this
theorem.

We firstly "lift the Theorem" on the level of complexes. We do it as
follows.

Let $A$ and $B$ be two associative algebras, and let $K$ be a dg
$B-A$-module (this means that it is a left $B$-module and left
$A$-module, and the left action commutes with the right action).
Define then a differential graded category with 2 objects, $a$ and
$b$, as follows. We set $\Mor(a,a)=A$, $\Mor(b,b)=B$, $\Mor(b,a)=K$,
$\Mor(a,b)=0$. To make this a dg category the only what we need is
that $A$ and $B$ are algebras, and $K$ is a $B-A$-module. Denote
this category by $\cat(A,B,K)$, see Section 5 for more detail.

Consider the Hochschild cohomological complex
$\Hoch^\mb(\cat(A,B,K))$ of this dg category. There are natural
projections $p_A\colon \Hoch^\mb(\cat(A,B,K))\to\Hoch^\mb(A)$ and
$p_B\colon\Hoch^\mb(\cat(A,B,K))\to\Hoch^\mb(B)$. The B.Keller's
theorem [Kel1] gives sufficient conditions for $p_A$ and $p_B$ being
quasi-isomorphisms. These conditions are that the natural maps
\begin{equation}\label{eqnewintro3}
B\rightarrow \mathrm{RHom}_{Mod-A}(K,K)
\end{equation}
and
\begin{equation}\label{eqnewintro4}
A^{opp}\rightarrow \mathrm{RHom}_{B-Mod}(K,K)
\end{equation}
are quasi-isomorphisms.

An easy computation shows that in the case when $A$ is Koszul
algebra, $B=A^{!opp}$ is the opposite to the Koszul dual algebra,
and $K$ is the Koszul complex of $A$, the Keller's condition is
satisfied (see Section 5).
\begin{remark}
According to the Remark at the end of Section 0.3, we should work with the $\mathbb{Z}$-graded category.
Therefore, our Hochschild complexes should be also compatible with this grading. More precisely, the Hochschild cochains should be finite sums of graded cochains. See Section 4.2.1 where it is explicitly stated.
\end{remark}
Consider the case when $A=S(V^*)\otimes\mathbb{C}[[\hbar]]$ and
$B=\Lambda(V)\otimes\mathbb{C}[[\hbar]]$. Denote in this case the
category $\cat(A,B,K)$ where $K$ is the Koszul complex of $A$, just
by $\cat$. Consider the following solid arrow diagram diagram:

\begin{equation}\label{eqnewintro5}
\xymatrix{&\Hoch^\mb(A)\\
T_\poly(V)\ar[ur]^{\mathcal{U}^S}\ar[dr]_{\mathcal{U}^\Lambda}\ar@{.>}[rr]^{\mathcal{F}}&&\Hoch^\mb(\cat)\ar[ul]_{p_A}\ar[dl]^{p_B}\\
&\Hoch^\mb(B)}
\end{equation}

The right "horn" was just defined. The maps $\U^S$ and $\U^\Lambda$
in the left "horn" are the following. We consider some {\it
universal} $L_\infty$ map $\U\colon T_\poly(V)\to\Hoch^\mb(S(V^*))$.
This means that we attribute some complex numbers to each Kontsevich
graph in his formality morphism in [K97], but which are not
necessarily the Kontsevich integrals (but the first Taylor
components is fixed, it is the Hochschild-Kostant-Rosenberg map).
The word "universal" again means that these numbers are the same for
all spaces $V$. Then we apply this $L_\infty$ morphism to our space
$V$, it is $\U^S$, and the composition of $\mathcal{D}\colon
T_\poly(V)\to T_\poly(V^*[1])$ with the $L_\infty$ morphism
$\U\colon T_\poly(V^*[1])\to \Hoch^\mb(\Lambda(V))$, constructed
from the same universal $L_\infty$ map.

The all solid arrows (\ref{eqnewintro5}) are quasiisomorphisms.
Therefore, they are homotopically invertible (see Section 4), and we
can speak about the homotopical commutativity of this diagram.

\begin{theorem*}
There exists a universal $L_\infty$ morphism $\U\colon
T_\poly(V)\to\Hoch^\mb(S(V^*))$ such that the solid arrow diagram
(\ref{eqnewintro5}) is homotopically commutative. The $L_\infty$
morphism corresponding by the Tamarkin's theory (see Sections 2 and 3) to
any choice of the Drinfeld associator satisfies this property.
\end{theorem*}

We firstly explain why our theorem about Koszul duality follows from
this theorem, and, secondly, how to prove this new theorem.

\subsection{}\label{section0.5}
If we know the homotopical commutativity of the solid arrow diagram
(of quasi-isomorphismes) (\ref{eqnewintro5}), we can construct the
dotted arrow $\mathcal{F}$ which is a $G_\infty$ quasi-isomorphism
$\F\colon T_\poly(V)\to\Hoch^\mb(\cat)$, which divides the diagram
into two homotopically commutative triangles. Then, if $\alpha$ is a
quadratic Poisson bivector field on $V$, the $L_\infty$ part of $\F$
defines a solution of the Maurer-Cartan equation $\F_*(\alpha)$ in
$\Hoch^\mb(\cat)$. A solution of the Maurer-Cartan equation in
$\Hoch^\mb(\cat)$ deforms the following four things simultaneously:
the algebra structures on $A$ and $B$, the differential in $K$, and
the bimodule structure. Using very simple arguments we then can
prove that this deformed complex $K$ is a free resolution of the
deformed $A$, and the deformed bimodule isomorphisms
(\ref{eqnewintro3})-(\ref{eqnewintro4}) give the Koszul duality
theorem. See Section 7 for detail.

\subsection{}\label{section0.6}
Here we outline the main ideas of Theorem \ref{section0.4}. First of
all, the two maps $p_A$ and $p_B$ in the right "horn" of the diagram
(\ref{eqnewintro5}) are maps of $B_\infty$ algebras (see [Kel1]).
Here $B_\infty$ is the braces dg operad, which acts on the
Hochschild cohomological complex of any algebra or dg category.
Formally it is defined as follows: a $B_\infty$ module structure on
$X$ is a dg bialgebra structure on the cofree coalgebra cogenerated
by $X[1]$ such that the coalgebra structure coincides with the given
one. The action of $B_\infty$ on the Hochschild complex
$\Hoch^\mb(A)$ of any dg algebra (or dg category) $A$ is constructed
by Getzler-Jones [GJ] via the braces operations.

Now define analogously the dg operad $B_\Lie$. A vector space $Y$ is
an algebra over $B_\Lie$ iff there is a dg Lie bialgebra structure
on the free Lie coalgebra cogenerated by $Y[1]$ such that the Lie
coalgebra structure coincides with the given one. The operads
$B_\Lie$ and $B_\infty$ are quasi-isomorphic by the Etingof-Kazhdan
(de)quantization. The construction of quasi-isomorphism of operads $B_\Lie\to
B_\infty$ depends on the choice of Drinfeld's associator [D].

The operad $B_\Lie$ is quasi-isomorphic to the Gerstenhaber
homotopical operad $G_\infty$, as is explained in [H], Section 6
(see also discussion in Section 2 of this paper). Finally, the
Gerstenhaber operad is Koszul, and $G_\infty$ is its Koszul
resolution. Thus, any $B_\infty$ algebra can be considered as a
$G_\infty$ algebra. As $G_\infty$ is a resolution of the
Gerstenhaber operad $G$, all three dg operads $G_\infty$,
$B_\infty$, and $B_\Lie$, are quasi-isomorphic to their cohomology
$G$.

(There is no canonical morphism from $B_\Lie$ to $B_\infty$. Any
such quasi-isomorphism gives a $G_\infty$ structure on the
Hochschild cohomological complex of any dg category. Any Drinfeld
associator [D] gives, via the Etingof-Kazhdan (de)quantization, such
a morphism of operads.)

Now consider the entire diagram (\ref{eqnewintro5}) as a diagram of
$G_\infty$ algebras and $G_\infty$ maps, where the $G_\infty$ action
on the Hochschild complexes is as above, it depends on the choice of
a map $B_\Lie\to B_\infty$. Then, if our diagram is homotopically
not commutative, it defines some $G_\infty$ automorphism of
$T_\poly(V)$.

This $G_\infty$ automorphism is clearly
$\mathrm{Aff}(V)$-equivariant. First of all, we prove that on the
level of cohomology the diagram (\ref{eqnewintro5}) is commutative.
This is in a sense the only new computation which we make in this
paper (see Section 5).

Thus, we can take the logarithm of this automorphism, which is
$G_\infty$-derivation. By the Tamarkin's $G_\infty$-ridigity of
$T_\poly(V)$, any $\mathrm{Aff}$-equivariant derivation is
homotopically inner. But any inner derivation acts non-trivially on
cohomology! On the other hand, a $G_\infty$-morphism homotopically
equivalent to identity, acts by the identity on cohomology. This
proves that our diagram is homotopically commutative. The only
property of this diagram which we have used is that it is defined
over $G_\infty$ and is commutative on the level of cohomology.

\subsection{}
When the author started to attack this problem, he started to prove
the homotopical commutativity of the diagram (\ref{eqnewintro5}) by
"physical" methods. Namely, the Kontsevich's formality in the
original proof given in [K97] is a particular case of the AKSZ model
on open disc [AKSZ], also called by Cattaneo and Felder the "Poisson
sigma-model". As usual in open theories, we should impose some
boundary conditions for the restrictions of the fields to the circle
$S^1=\partial D^2$. Maxim Kontsevich considers the boundary
condition $"p=0"$ on all arcs. This, certainly together with other
mathematical insights, led him in [K97] to the formality theorem.

Our idea was to divide $S^1$ by two parts, fixing two points $\{0\}$
and $\{\infty\}$ (in the Kontsevich's case only $\{\infty\}$ is
fixed). Then, we impose the boundary condition $"x=0"$ on all left
arcs, and $"p=0"$ on all right arcs. This seems to be very
reasonable, and the author hoped to construct in this way an
$L_\infty$ quasi-isomorphism $\F$ (the dotted arrow in
(\ref{eqnewintro5})), making the two triangles homotopically
commutative.

Somehow, he did not succeed in that. From the point of view of this
paper, it seems that the reason for that is the following.

The author worked with the Kontsevich's propagator in [K97], namely,
with

\begin{equation}\label{eqpropend1}
\varphi(z_1,z_2)=\frac1{2\pi
i}d\Log\frac{(z_1-z_2)(z_1-\overline{z}_2)}{(\overline{z}_1-z_2)(\overline{z}_1-\overline{z}_2)}
\end{equation}
(here $z_1$ and $z_2$ are distinct points of the complex upper
half-plane).

In this paper we deal with the Tamarkin's quantization.
Conjecturally (see [K99]) when this formality morphism is
constructed from the Knizhnik-Zamolodchikov Drinfeld's associator,
it coincides (as a universal $L_\infty$ morphism, see above) with
the $L_\infty$ morphism, constructed from the "another Kontsevich's
propagator". This is "the half" of (\ref{eqpropend1}):

\begin{equation}\label{eqpropend2}
\varphi_1(z_1,z_2)=\frac1{2\pi
i}d\Log\frac{z_1-z_2}{\overline{z}_1-z_2}
\end{equation}

Kontsevich proved (unpublished) that this propagator also leads to
an $L_\infty$ morphism from $T_\poly(V)$ to $\Hoch^\mb(S(V^*))$. If
this conjecture (that the Tamarkin's theory in the
Knizhnik-Zamolodchikov case gives this propagator) is true, we
should try to elaborate the physical idea described above (with the
two boundary conditions) for this propagator. The reason is that it
is not a priori clear that the Kontsevich's first propagator
$\varphi(z_1,z_2)$ comes from any Drinfeld's associator, and
therefore from the Tamarkin's theory.

We are going to come back to these questions in the sequel.

\subsection{}
We tried to make the exposition as self-contained as possible. In
particular, we prove in Section 2.4 the main Lemma in the Tamarkin's
proof [T1] of the Kontsevich formality, because we use it here not
only for the first cohomology as in [T1] and [H], and also for 0-th
cohomology. We give a simple proof of it for all cohomology for
completeness. As well, we reproduce in Section 4.2 the proof of
Keller's theorem from [Kel1], because in [Kel1] some details are
omitted. Nevertheless, in one point we did not overcome some
vagueness. This is the using of the homotopical relation for maps of
dg operads or algebras over dg operads. Some implications like
"homotopical maps of dg operads induce homotopical morphisms of
algebras" in Section 3, are stated without proofs. Finally in
Section 5 we give a construction of the homotopical category of dg
Lie algebras through the "right cylinder" in the sense of [Q], which
is suitable for the proof of the Main Theorem in Section 7.

\subsection{}
The paper is organized as follows:

In Section 1 we develop the Koszul duality for algebras over a
discrete valuation rings. Our main example is the algebras over
$\mathbb{C}[[\hbar]]$, and we should justify that the main theorems
of Koszul duality for associative algebras hold in this context;

In Section 2 we give a brief exposition of the Tamarkin's theory
[T1]. The Hinich's paper [H] is a very good survey, but we achieve
some more clarity in the computation of deformation cohomology of
$T_\poly(V)$ over the operad $G_\infty$ of homotopy Gerstenhaber
algebras. As well, in the Appendix after Section 2.5 we give a
deduction of the existence of Kontsevich formality over $\mathbb{Q}$
from its existence over $\mathbb{C}$, which differs from the
Drinfeld's approach [D]. This deduction seems to be new;

In Section 3 we touch some unsolved problems in the Tamarkin's
theory and leave them unsolved, wee only need to know here that any
map of operads $G_\infty\to B_\infty$ defined up to homotopy,
defines a {\it universal} $G_\infty$ map $T_\poly(V)\to
\Hoch^\mb(S(V^*))$ where the $G_\infty$ structure on
$\Hoch^\mb(S(V^*))$ is defined via the map of operads. The
deformation quantizations for which our Main Theorem is true belong
to the image of the map $\mathfrak{X}$ defined there;

In Section 4 for introduce differential graded categories, give a
construction of the Keller's dg category from [Kel1] associated with
a Keller's triple, and reformulate our Main Theorem in this new
setting. We get a more general statement, which is, however, more
easy to prove;

A very short Section 5 is just a place to relax before the long
computation in Section 6, here we recall the explicit construction
[Sh3] of the Quillen's homotopical category via the right cylinder.
The advantage of this construction is that it is immediately clear
from it that two homotopical $L_\infty$ maps map a solution of the
Maurer-Cartan equation to gauge equivalent solutions (Lemma 5.2);

In Section 6 we construct the Hochschild-Kostant-Rosenberg map for
the Keller's dg category. This computation is done in terms of
graphs, closed to the ones from [K97]. Originally the author got
this computation truing to construct the $L_\infty$ morphism
$\F\colon T_\poly(V)\to\Hoch^\mb(\cat(A,B,K))$ dividing the diagram
(\ref{eqnewintro5}) into two commutative triangles, by "physical"
methods. The computation here is the only what the author succeed to
do in this direction;

The final Section 7 is the culmination of our story. Here we deduce
the Main Theorem on Koszul duality in deformation quantization from
Theorem 4.4. The idea is that from the diagram (\ref{eqnewintro5})
associates with a quadratic Poisson bivector $\alpha$ on $V$ a
solution of the Maurer-Cartan equation in the Hochschild complex of
the Keller's dg category. This Maurer-Cartan elements defines an
$A_\infty$ deformation of the Keller's category, and, in particular,
deforms the Koszul complex. This is enough to conclude that the two
deformed algebras are Koszul dual.

\subsubsection*{Acknowledgements}
I am very grateful to Maxim Kontsevich who taught me his formality
theorem and many related topics. Discussions with Victor Ginzburg,
Pavel Etingof and Bernhard Keller were very valuable for me. It was
Victor Ginzburg who put my attention on the assumptions on $A_0$ in
the theory of Koszul duality, related with flatness of
$A_0$-modules, and explained to me why $A_0$ is semisimple in [BGS].
And it was Bernhard Keller who explained to me in our correspondence
some foundations about dg categories, as well as his constructions
from [Kel1]. But more than to the others, I am indebted to Dima
Tamarkin. Discussions with Dima after my talk at the Nothwestern
University on the subject of the paper, and thereafter, shed new
light to many of my previous constructions, and finally helped me to
prove the Main Theorem of this paper.

I express my gratitude to the MIT and to the University of Chicago
which I visited in October-November 2007 and where a part of this
work was done, for a very stimulating atmosphere and for the
possibility of many valuable discussions, as well as for their
hospitality and particular financial support.

I am grateful to the research grant R1F105L15 of the University of
Luxembourg for partial financial support.

\section{Koszul duality for algebras over a discrete valuation
ring}\label{section1} Here we give a brief overview of the Koszul
duality. Our main reference is Section 2 of [BGS]. In loc.cit., the
zero degree component $A_0$ is supposed to be a (non-commutative)
semisimple algebra over the base field $k$. For our applications in
deformation quantization, we should consider algebras over
$\mathbb{C}[[\hbar]]$. For this reason, we show that the theory of
Koszul duality may be defined over an arbitrary commutative discrete
valuation ring. This result seems to be new, although L.Positselski
announced in [P] that the zero degree component $A_0$ may be an
arbitrary algebra over the base field.
\subsection{}\label{section1.1}
The main classical example of Koszul dual algebras are the algebras
$A=S(V^*)$ and $B=\Lambda(V)$, where $V$ is a finite-dimensional
vector space over the base field $k$. In general, suppose $A_0$ is a
fixed $k$-algebra. Koszulness is a property of a graded algebra
\begin{equation}\label{01.1}
A=A_0\oplus A_1\oplus A_2\oplus\dots
\end{equation}
that is,
\begin{equation}\label{01.2}
A_i\cdot A_j\subset A_{i+j}
\end{equation}
In our example with $S(V^*)$ and $\Lambda(V)$ the algebra $A_0=k$,
it is the simplest possible case. In general, all $A_i$ are
$A_0$-bimodules.

There is a natural projection $p\colon A\to A_0$ which endows $A_0$
with a (left) $A$-module structure. Denote by $A-\Mod$ the category
of all left $A$-modules, and by $A-\mathrm{mod}$ the category of
graded left $A$-modules.

The $A$-module $A_0$ always has a free resolution in
$A-\mathrm{mod}$
\begin{equation}\label{01.3}
\dots\rightarrow M_2\rightarrow M_1\rightarrow M_0\rightarrow 0
\end{equation}
such that $M_i$ is a graded $A$-module generated by elements of
degrees $\ge i$. Indeed, the bar-resolution
\begin{equation}\label{01.4}
\dots \rightarrow A\otimes_k A_+^{\otimes 2}\otimes_k A_0\rightarrow
A\otimes_k A_+\otimes_k A_0\rightarrow A\otimes_k A_0\rightarrow 0
\end{equation}
obeys this property. (Here $A_+=A_1\oplus A_2\oplus\dots$). This
motivates the following definition:
\begin{definition}
A graded algebra (\ref{01.1}) is called {\it Koszul} if the
$A$-module $A_0$ admits a projective resolution (\ref{01.3}) in
$A-\mathrm{mod}$ such that each $M_i$ is finitely generated by
elements of degree $i$.
\end{definition}
For our example with the symmetric and the exterior algebra, such a
resolution is the following {\it Koszul complex}:
\begin{equation}\label{01.5}
\dots\rightarrow S(V^*)\otimes\Lambda^3(V)^*\rightarrow
S(V^*)\otimes\Lambda^2(V)^*\rightarrow S(V^*)\otimes V^*\rightarrow
S(V^*)\rightarrow 0
\end{equation}
with the differential
\begin{equation}\label{01.6}
d=\sum_{i=1}^{\dim V}x_i\otimes \frac{\partial}{\partial\xi_i}
\end{equation}
Here $\{x_i\}$ is a basis in the vector space $V^*$ and $\{\xi_i\}$
is the corresponding basis in $V^*[1]$. The differential is a
$\gl(V)$-invariant element, it does not depend on the choice of
basis $\{x_i\}$ of the vector space $V^*$.

\subsection{}\label{section1.2}
Here we explain some consequences of the definition of Koszul
algebra, leading to the concept of Koszul duality for quadratic
algebras. In Sections \ref{section1.2.1}-\ref{section1.2.3} $A_0$
may be arbitrary finite-dimensional algebra over the ground field
$k$, and in Sections 1.2.4-1.2.6 we suppose that  $A_0$ is a
semisimple finite-dimensional algebra over $k$ (see [BGS]).

\subsubsection{}\label{section1.2.1}
Let $A$ be a graded algebra. Then the space
$\Ext^\mb_{A-\Mod}(A_0,A_0)$ is naturally {\it bigraded}. We write
$\Ext^n_{A-\Mod}(A_0,A_0)=\oplus_{a+b=n}\Ext^{a,b}(A_0,A_0)$. From
the bar-resolution (\ref{01.4}) we see that for a general algebra
$A$ the only non-zero $\Ext^{a,b}(A_0,A_0)$ appear for $a\le -b$
(here $a$ is the cohomological grading and $b$ is the inner
grading). In the Koszul case the only nonzero summands are
$\Ext^{a,-a}$. Let us analyze this condition for $a=1$ and $a=2$.

\subsubsection{}\label{section1.2.2}
\begin{lemma}
Suppose $A$ is a graded algebra.
\begin{itemize}
\item[1.] If $\Ext_{A-\Mod}^1(A_0,A_0)=\Ext^{1,-1}(A_0,A_0)$ (that is, all
$\Ext^{1,-b}=0$ for $b>1$), the algebra $A$ is 1-generated. The
latter means that the algebra $A$ in the form of (\ref{01.1}) is
generated over $A_0$ by $A_1$;
\item[2.] if, furthermore,
$\Ext_{A-\Mod}^2(A_0,A_0)=\Ext^{2,-2}(A_0,A_0)$ (that is,
$\Ext^{2,-\ell}(A_0,A_0)=0$ for $\ell\ge 3$), the algebra $A$ is
quadratic. This means that $A=T_{A_0}(A_1)/I$ where $I$ is a graded
ideal generated as a two-sided ideal by $I_2=I\cap A_2$.
\end{itemize}
\end{lemma}
See [BGS], Section 2.3.

This Lemma implies that any Koszul algebra is quadratic. So, in fact
the Koszulness is a property of quadratic algebras.

\subsubsection{}\label{section1.2.3}
From now on, we use the notation $I=I_2$ for the intersection of the
graded ideal $I$ in $T_{A_0}(A_1)$ with $A_2$. Any quadratic algebra
is uniquely defined by the triple $(A_0,A_1,I\subset
A_1\otimes_{A_0}A_1)$.

Using the bar-complex (\ref{01.4}), it is very easy to compute the
"diagonal part" $\oplus_\ell\Ext^{\ell,-\ell}(A_0,A_0)\subset
\Ext^\mb_{A-\Mod}(A_0,A_0)$ for any algebra $A$. Let us formulate
the answer.

Define from a triple $(A_0, A_1, I)$ another triple $(A_0^\vee,
A_1^\vee, I^\vee)$, as follows. Suppose $A_1$ and $I$ are flat
$A_0$-bimodules. We set $A_0^\vee=A_0$, $A_1^\vee=\Hom_{A_0}(A_1,
A_0)[-1]$. Define now $I^\vee$. Denote firstly
$A_1^*=\Hom_{A_0}(A_1,A_0)$. There is a pairing
$(A_1\otimes_{A_0}A_1)\otimes (A_1^*\otimes_{A_0} A_1^*)\to A_0$ which
is non-degenerate. Denote by $I^*$ the subspace in
$A_1^*\otimes_{A_0}A_1^*$ dual to $I$. Denote by $I^\vee=I^*[-2]$,
it is a subspace in $A_1^\vee\otimes_{A_0}A_1^\vee$. The triple
$(A_0, A_1^\vee, I^\vee)$ generates some quadratic algebra, denote
it by $A^\vee$.

Let now $A$ be any 1-generated not necessarily quadratic algebra.
Then the {\it quadratic part} $A^q$ is well-defined. Let $A$ be a
quotient of $T_{A_0}(A_1)$ by graded not necessarily quadratic
ideal. We define $A^q$ as the quadratic algebra associated with the
triple $(A_0, A_1, I\cap A_2)$. There is a canonical surjection
$A^q\to A$ which is an isomorphism in degrees 0, 1, and 2.

\begin{lemma}
Let $A$ be a 1-generated algebra over $A_0$. Then the diagonal
cohomology $\oplus_{\ell}\Ext^{\ell,-\ell}(A_0,A_0)$ as algebra is
canonically isomorphic to the algebra opposed to $(A^q)^\vee$. Here
by the opposed algebra to an algebra $B$ we understand the product
$b_1\star^{\opp}b_2=b_2\star b_1$.
\end{lemma}

It is a direct consequence from the bar-resolution (\ref{01.4}).

In particular, let now a graded 1-generated algebra $A$ be Koszul.
Then $\Ext^\mb_{A-\Mod}(A_0,A_0)=(A^\vee)^{\opp}$. This follows from
the identity $A=A^q$ for a quadratic algebra $A$ (in particular, for
Koszul $A$), and from the equality of the all $\Ext$s to its
diagonal part for any Koszul algebra.

\subsubsection{}\label{section1.2.4}
The inverse is also true, under an assumption on $A_0$.
\begin{lemma}
Suppose $A_0$ is a simple finite-dimensional algebra over $k$ and
$A$ is a quadratic algebra over $A_0$. Then if
$\Ext_{A-\Mod}^\mb(A_0,A_0)$ is equal to its diagonal part, then $A$
is Koszul. In particular, if
$\Ext^\mb_{A-\Mod}(A_0,A_0)=(A^{q\vee})^{\opp}$, then $A$ is Koszul.
\end{lemma}

See [BGS], Proposition 2.1.3.

Let us comment why we need here a condition on $A_0$. The projective
resolution which we need to prove that $A$ is Koszul is constructed
inductively. We construct a resolution
\begin{equation}\label{01.10}
\dots\rightarrow P_3\rightarrow P_2\rightarrow P_1\rightarrow
P_0\rightarrow 0
\end{equation}
satisfying the property of Definition \ref{section1.2} and such that
the differential is injective on $P_i^i$. We set $P_0=A$. To perform
the step of induction, set $K=\mathrm{ker}(P_i\to P_{i-1})$. We
have: $\Ext_{A-{\mathrm{Mod}}}^{i+1}(A_0, A_0)=\Hom_{A-\Mod}(K,
A_0)$. From the condition of lemma we conclude that $K$ is generated
by the part $K^{i+1}$ of inner degree $i+1$ (here for simplicity we
suppose that $A$ has trivial cohomological grading). Then we put
$P_{i+1}=A\otimes_{A_0} K^{i+1}$. But then we need to check that the
image  of the map $P_{i+1}\to P_i$ is $K$. For this we necessarily
need to know that $K^{i+1}$ is a flat left $A_0$-module. For this it
is sufficiently to know that $K$ is. So we need a theorem like the
following: the kernel of any map of good (flat, etc.) $A_0$-modules
is again a flat $A_0$-module. It does not follow from any general
things, it is a property of $A_0$. It is the case if {\it any}
module is flat, as in the case of a finite-dimensional simple
algebra. For another possible condition, see Section
\ref{section1.3}.

\subsubsection{}\label{section1.2.5}
\begin{proposition}
Suppose $A$ is a quadratic algebra defined from a triple
$(A_0,A_1,I)$ where $A_1$ and $I$ are flat $A_0$-bimodules. Suppose
$A$ is Koszul. Then $(A^\vee)^\opp$ is also Koszul.
\end{proposition}

\begin{remark}
It is clear that $A$ is Koszul iff $A^\opp$ is Koszul.
\end{remark}

We give a sketch of proof, which is essentially given by the
construction of the Koszul complex. For $A=S(V^*)$ the Koszul
complex is constructed in (\ref{01.5}).

Let $A=(A_0,A_1,I)$ be a quadratic algebra, and let $A_1$ and
$I\subset A_1\otimes_{A_0}A_1$ be flat $A_0$-bimodules. We define
the Koszul complex
\begin{equation}\label{01.20}
\dots \rightarrow K_3\rightarrow K_2\rightarrow K_1\rightarrow
K_0\rightarrow 0
\end{equation}
We set
\begin{equation}\label{01.21}
K_i=A\otimes_{A_0} K_i^i
\end{equation}
where
\begin{equation}\label{01.22}
K_i^i=\bigcap_\ell A_1^{\otimes\ell}\otimes_{A_0} I\otimes_{A_0}
A_1^{\otimes i-\ell-2}
\end{equation}
In particular, $K_0^0=A_0$, $K_1^1=A_1$, $K_2^2=I$. The differential
$d\colon K_i\to K_{i-1}$ is defined as the restriction of the map
$\hat{d}\colon A\otimes_{A_0} A_1^{\otimes i}\to
A\otimes_{A_0}A_1^{\otimes (i-1)}$ given as
\begin{equation}\label{01.23}
a\otimes v_1\otimes\dots\otimes v_i\mapsto (av_1)\otimes
v_2\otimes\dots\otimes v_{i-1}
\end{equation}
Clearly $d^2=0$. The complex (\ref{01.20}) is called the {\it Koszul
complex} of the quadratic algebra $A$.

\begin{lemma}
Let $A=(A_0,A_1,I)$ be a quadratic algebra, $A_1$ and $I$ flat
$A_0$-bimodules. Suppose, additionally, that $A_0$ is a
finite-dimensional semisimple algebra over $k$. Then its Koszul
complex is acyclic except degree 0 iff $A$ is Koszul.
\end{lemma}
See [BGS], Theorem 2.6.1 for a proof. {\it In the proof it is
essential that the modules $K_i^i$ a flat left $A_0$-bimodules}. In
general the tensor product of two flat modules is flat, but there is
no theorem which guarantees the same about the intersection of two
flat submodules. In the case which is considered in [BGS], any
module over a finite-dimensional semisimple algebra is flat.

Let us note that the part "only if" also follows from Lemma
\ref{section1.2.4}.

The Proposition follows from this Lemma easily.

Indeed, it is clear that $K_i=A\otimes_{A_0} [(A^!)^*]^i[-i]$ and
that the Koszul complex $K$ of a Koszul algebra satisfies the
Definition \ref{section1.1}. Then the dual complex $K^*$ also
satisfies the Definition \ref{section1.1} and can be written as
$K^*=A^*\otimes_{A_0} A^!$. We immediately check that it coincides
with the Koszul complex of the quadratic algebra $A^!$ because
$(A^!)^!=A$ for any quadratic algebra $A$. Then from its acyclicity
follows that $A^!$ is Koszul.

\subsubsection{}\label{section1.2.6}
We summarize the discussion above in the following theorem.
\begin{theorem*}
Let $A=(A_0,A_1,I)$ be a quadratic algebra, $A_1$ and $I$ be flat
$A_0$-bimodules, and $A_0$ be semisimple finite-dimensional algebra
over $k$. Then $A$ is Koszul if and only if the quadratic dual $A^!$
is also Koszul, and in this case
\begin{equation}\label{01.25}
\Ext^i_{A-\Mod}(A_0,A_0)=[(A^!)^{\opp}]_i[-i]
\end{equation}
and
\begin{equation}\label{01.26}
\Ext_{A^!-\Mod}^i(A_0,A_0)=[A^\opp]_i[-i]
\end{equation}
for any integer $i\ge 0$. \qed
\end{theorem*}

\subsection{}\label{section1.3}
In the context of deformation quantization, all our algebras are
algebras over the formal power series $\mathbb{C}[[\hbar]]$,
therefore, $A_0=\mathbb{C}[[\hbar]]$. The theory of Koszul algebras
as it is developed in [BGS] does not cover this case. In this
Subsection we explain that in the theory of Koszul duality $A_0$ may
be any commutative discrete valuation ring (see [AM], [M]).
L.Positselski announced in [P] that $A_0$ may be any algebra over
$k$.

Recall the definition of a discrete valuation ring.

\begin{definition}
A commutative ring is called a discrete valuation ring if it is an
integrally closed domain with only one nonzero prime ideal. In
particular, a discrete valuation ring is a local ring.
\end{definition}

The two main examples are the following:

(1) Let $C$ be an affine algebraic curve, and let $p\in C$ be a
non-singular point (not necessarily closed). Then the local ring
$\mathcal{O}_p(C)$ is a discrete valuation ring (recall that in
dimension 1 integrally closed=nonsingular);

(2) let $C$ and $p$ be as above; we can consider the completion of
the local ring $\mathcal{O}_p(C)$ by the powers of the maximal
ideal. Denote this ring by $\widehat{\mathcal{O}}_p(C)$, this is a
discrete valuation ring. In particular, $k[[\hbar]]$ is a discrete
valuation ring.

It is known that any discrete valuation ring is Noetherian and is a
principal ideal domain (see [M], Theorem 11.1).

To extend the theory of Section \ref{section1.2} to the case when
$A_0$ is a discrete valuation ring we need to know that the
intersection of flat submodules  over a discrete valuation ring
(Section \ref{section1.2.5}), and the kernel of a map of flat
modules over a discrete valuation ring (Section \ref{section1.2.4})
are flat. This is guaranteed by the following, more general, result:

\begin{lemma}
Let $R$ be a discrete valuation ring, and let $M$ be a flat
$R$-module. Then any submodule of $M$ is again flat.
\begin{proof}
Let $R$ be a ring and $M$ is an $R$-module. Then $M$ is flat if and
only if for any finitely generated ideal $I\subset R$ the natural
map $I\otimes_R M\to R\otimes_R M$ is injective (see [M], Theorem
7.7). Any ideal in a discrete valuation ring is principal ([M],
Theorem 11.1); therefore flatness of a module over a discrete
valuation ring is the same that torsion-free (a module $M$ is called
torsion-free if $x\ne 0$, $m\ne 0$ implies $xm\ne 0$). So now our
Lemma follows from the fact that a submodule over a torsion-free
module is torsion-free.
\end{proof}
\end{lemma}

\begin{remark}
If $R$ is any local ring and $M$ is a finite $R$-module, then
flatness of $M$ implies that $M$ is free ([M], Theorem 7.10).
Nevertheless, in dimension $\ge 2$ a submodule of a free module may
be not free. For example, one can take the (localization of the)
coordinate ring of a curve in $\mathbb{A}^2$.
\end{remark}

Combining the Lemma above with the discussion of Section
\ref{section1.2}, we get the following Theorem:
\begin{theorem*}
Let $A=(A_0,A_1,I)$ be a quadratic algebra, with $A_0$ a commutative
discrete valuation ring, and $A_1$, $I$ flat $A_0$-modules. Then $A$
is Koszul iff $A^!$ is, and in this case
\begin{equation}\label{01.27}
\Ext^i_{A-\Mod}(A_0,A_0)=[(A^!)^{\opp}]_i[-i]
\end{equation}
and
\begin{equation}\label{01.28}
\Ext_{A^!-\Mod}^i(A_0,A_0)=[A^\opp]_i[-i]
\end{equation}
for any integer $i\ge 0$. \qed
\end{theorem*}
We will use this Theorem only for $A_0=k[[\hbar]]$.

\begin{remark}
This theorem has an analogue for Dedekind domains. Namely, the
localization of a Dedekind domain at any prime ideal is a discrete
valuation ring, this is a global version of it. The main example of
a Dedekind domain is the coordinate ring of a non-singular affine
curve. Suppose $A_0$ is a Dedekind domain. We say that a quadratic
algebra over $A_0$ is Koszul if its localization at any prime ideal
is Koszul. Then we can prove the theorem analogous to the above for
$A_0$ a Dedekind domain. More generally, we can speak about sheaves
of Koszul dual quadratic algebras. At the moment the author does not
know any interesting example of such situation, but he does not
doubt that these examples exist.
\end{remark}

\begin{remark}
Leonid Positselski claims in [P] that he constructed the analogous
theory for any $A_0$. The arguments in [P] are rather complicated
comparably with ours'; for the readers' convenience, we gave here a
more direct simple proof in the case of discrete valuation rings.
\end{remark}

\section{Tamarkin's approach to the Kontsevich
formality}\label{section2} Here we overview the Tamarkin's proof of
Kontsevich formality theorem. The main references are [T1] and [H],
some variations which allow to avoid the using of the
Etingof-Kazhdan quantization (but replace it by another
transcendental construction) were made by Kontsevich in [K99].

\subsection{Kontsevich formality}\label{section2.1}
For any associative algebra $A$ we denote by $\Hoch^\mb(A)$ the
cohomological Hochschild complex of $A$. When $A=C^\infty(M)$ is the
algebra of smooth functions on a smooth manifold $M$, we consider
some completed tensor powers, or the polydifferential part of the
usual Hochschild complex (see, e.g., [K97]). Under this assumption,
the Hochschild cohomology of $A=C^\infty(M)$ is equal to smooth
polyvector fields $T_\poly(M)$. More precisely, consider the
following Hochschild-Kostant-Rosenberg map $\varphi\colon
T_\poly(M)\to\Hoch^\mb(C^\infty(M))$:
\begin{equation}\label{2.1.1}
\varphi(\gamma)=\{f_1\otimes\dots\otimes f_k\mapsto
\frac1{k!}\gamma(df_1\wedge\dots\wedge df_k)\}
\end{equation}
for $\gamma$ a $k$-polyvector field. Then the
Hochschild-Kostant-Rosenberg theorem is
\begin{lemma}
\begin{itemize}
\item[1.] For any polyvector field $\gamma$, the cochain
$\varphi(\gamma)$ is a cocycle; this gives an isomorphism of
(completed or polydifferential) Hochschild cohomology of
$A=C^\infty(M)$ with $T_\poly(M)$;
\item[2.] the bracket induced on the Hochschild cohomology from the
Gerstenhaber bracket coincides, via the map $\varphi$, with the
Schouten-Nijenhuis bracket of polyvector fields.
\end{itemize}
\end{lemma}
See, e.g., [K97] for definition of the Gerstenhaber and
Schouten-Nijenhuis brackets.

The second claim of the Lemma means that
\begin{equation}\label{2.1.2}
[\varphi(\gamma_1),\varphi(\gamma_2)]_G=\varphi([\gamma_1,\gamma_2]_{SN})+d_\Hoch\U_2(\gamma_1,\gamma_2)
\end{equation}
for some $\U_2\colon
\Lambda^2(T_\poly(M))\to\Hoch^\mb(C^\infty(M))[-1]$ (we denoted by
$[\,\ ]_G$ the Gerstenhaber bracket and by $[\ ,\ ]_{SN}$ the
Schouten-Nijenhuis bracket).

In the case when $M=\mathbb{C}^d$ M.Kontsevich constructed in [K97]
an $L_\infty$ morphism $\U\colon T_\poly(\mathbb{C}^d)\to
\Hoch^\mb(S(\mathbb{C}^{d*}))$ whose first Taylor component is the
Hochschild-Kostant-Rosenberg map $\varphi$. (Here we consider {\it
polynomial} polyvector fields, and there is no necessity to complete
the Hochschild complex). The second Taylor component $\U_2$ should
then satisfy (\ref{2.1.2}), and so on. This result is called {\it
the Kontsevich's formality theorem}.

(The result for a general manifold $M$ can be deduced from this
local statement, see [K97], Section 7).

The original Kontsevich's proof uses ideas of topological field
theory, namely, the Alexandrov-Kontsevich-Schwarz-Zaboronsky (AKSZ)
model, see [AKSZ]. Therefore, some transcendental complex numbers,
the "Feynmann integrals" of the theory, are involved into the
construction. The Kontsevich's proof appeared in 1997.

One year later, in 1998, D.Tamarkin found in [T] another proof of
the Kontsevich formality for $\mathbb{C}^d$, using absolutely
different technique. In the rest of this Section we outline the
Tamarkin's proof [T], [H] in the form we use it in the sequel.

\subsection{The idea of the Tamarkin's proof}\label{section2.2}
The main idea it to construct not only an $L_\infty$ map from
$T_\poly(\mathbb{C}^d)$ to $\Hoch^\mb(S(\mathbb{C}^{d*}))$ but to
involve the entire structure on polyvector fields and the Hochschild
complex. This is the structure of (homotopy) Gerstenhaber algebra.
For example, on polyvector fields (on any manifold) one has two
operations: the wedge product $\gamma_1\wedge \gamma_1$ of degree 0,
and the Lie bracket $[\gamma_1,\gamma_2]_{SN}$ of degree $-1$, and
they are compatible as
\begin{equation}\label{2.2.1}
[\gamma_1,\gamma_2\wedge\gamma_3]=[\gamma_1,\gamma_2]\wedge\gamma_2\pm
\gamma_2\wedge[\gamma_1,\gamma_3]
\end{equation}
This is called a Gerstenhaber algebra. To consider
$T_\poly(\mathbb{C}^d)$ as a Gerstenhaber algebra simplifies the
problem because of the following Lemma:
\begin{lemma}
The polyvector fields $T_\poly(\mathbb{C}^d)$ is rigid as a homotopy
Gerstenhaber algebra. More precisely, any
$\mathrm{Aff}(\mathbb{C}^d)$-equivariant deformation of
$T_\poly(\mathbb{C}^d)$ as a homotopy Gerstenhaber algebra is
homotopically equivalent to trivial deformation.
\end{lemma}

We should explain what these words mean, we do it in the next
Subsections. Let us now explain how it helps to prove the
Kontsevich's formality theorem.

It is true, and technically it is the hardest place in the proof,
that there is a structure of homotopical Gerstenhaber algebra (see
Section \ref{section2.3}) on the Hochschild complex $\Hoch^\mb(A)$
of any associative algebra $A$. It is non-trivial, because the
cup-product of Hochschild cochains $\Psi_1\cup\Psi_2$ and the
Gerstenhaber bracket $[\Psi_1,\Psi_2]_G$ do not obey the
compatibility (\ref{2.2.1}), as it should be in a Gerstenhaber
algebra. It obeys it only up to a homotopy, and to find explicitly
this structure uses also either some integrals like in [K99], or
Drinfeld's Knizhnik-Zamolodchikov iterated integrals, as in [T]. We
discuss it in Section \ref{section2.5}. Now suppose that this
structure exists, such that the Lie and Commutative parts of this
structure are equivalent to the Gerstenhaber bracket and the
cup-product on Hochschild cochains.

Then, as usual in homotopical algebra, there exists a homotopical
Gerstenhaber algebra structure on the cohomology, equivalent to this
structure on the cochains (it is something like "Massey operations"
by Merkulov and Kontsevich-Soibelman). This push-forwarded structure
is uniquely defined up to homotopy.

Now we can consider this structure as a formal deformation of the
classical pure Gerstenhaber algebra on $T_\poly(\mathbb{C}^d)$.
Indeed, we rescale the Taylor components of this structure, such
that the weight of $k$-linear Taylor components is $\lambda^{k-2}$.
This gives again a homotopical Gerstenhaber structure, which value
at $\lambda=0$ is the classical Gerstenhaber structure on polyvector
fields, because of the compatibility with Lie and Commutative
structure, and by the Hochschild-Kostant-Rosenberg theorem.

Now we apply the Lemma above. All steps of our construction are
$\mathrm{Aff}(\mathbb{C}^d)$-invariant, therefore, the obtained
deformation can be chosen $\mathrm{Aff}(\mathbb{C}^d)$-equivariant.
Then the Lemma says that this deformation is trivial, and the two
homotopical Gerstenhaber structures on $T_\poly(\mathbb{C}^d)$ are
in fact isomorphic. This implies the Kontsevich's formality in the
stronger, Gerstenhaber algebra isomorphism, form.

\subsection{Koszul operads}\label{section2.3}
From our point of view, the Koszulness of an operad $\mathcal{P}$ is
very important because in this case any $\mathcal{P}$-algebra $A$
has "very economic" resolution which is free dg
$\mathcal{P}$-algebra. In the case of the operad
$\mathcal{P}=Assoc$, this "very economic" resolution is the
Quillen's bar-cobar construction. Thereafter, we use this free
resolution to compute the (truncated) deformation complex of $A$ as
$\mathcal{P}$-algebra. In the case of $\mathcal{P}=Assoc$ this
deformation complex is the Hochschild cohomological complex of $A$
without the zero degree term, that is $\Hoch^\mb(A)/A$.

We will consider only operads of dg $\mathbb{C}$-vector spaces here,
with one of the two possible symmetric monoidal structures. A {\it
quadratic} operad generated by a vector space $E$ over $\mathbb{C}$
with an action of the symmetric group $\Sigma_2$ of two variables,
with a $\Sigma_3$-invariant space of relations $R\subset
\mathrm{Ind}_{\Sigma_2}^{\Sigma_3}E\otimes E$ (here $\Sigma_2$ acts
only on the second factor) is the quotient of the free operad
$\mathcal{P}$ generated by $\mathcal{P}(2)=E$ by the space of
relations $R\subset \mathcal{P}(3)$. The operads $Lie$, $Comm$,
$Assoc$ are quadratic, as well as the Gerstenhaber and the Poisson
operads. See [GK], Section 2.1 for more detail. For a quadratic
operad $\mathcal{P}$ define the quadratic dual operad
$\mathcal{P}^!$ as the quadratic operad generated by
$\mathcal{P}^!(2)=E^*[1]$, with the space of relations $R^*$ in
$\mathrm{Ind}_{\Sigma_2}^{\Sigma_3}E^*[1]\otimes E^*[1]$ equal to
the orthogonal complement to
$R\subset\mathrm{Ind}_{\Sigma_2}^{\Sigma_3}E\otimes E$. Example:
$Com^!=Lie[-1]$, $Assoc^!=Assoc[-1]$,
$(\mathcal{P}^!)^!=\mathcal{P}$.

Let $\mathcal{P}$ be a general, not necessarily quadratic, operad.
For simplicity, we suppose that all vector spaces $\mathcal{P}(n)$
of an operad $\mathcal{P}$ are finite-dimensional. Recall the
construction of the bar-complex of $\mathcal{P}$, see [GK], Section
3.2. Denote the bar complex of $\mathcal{P}$ by
$\mathbf{D}(\mathcal{P})$. Then one has:
$\mathbf{D}(\mathbf{D}(\mathcal{P}))$ is quasi-isomorphic to
$\mathcal{P}$ ([GK], Theorem 3.2.16). Let now $\mathcal{P}$ be a
quadratic operad. Then the bar-complex $\mathbf{D}(\mathcal{P})$ is
a negatively-graded dg operad whose 0-th cohomology is canonically
the quadratic dual operad $\mathcal{P}^!$. A quadratic operad
$\mathcal{P}$ is called {\it Koszul} if the bar-complex
$\mathbf{D}(\mathcal{P})$ is a resolution of $\mathcal{P}^!$. In
this case $\mathbf{D}(\mathcal{P}^!)$ gives a free resolution of the
operad $\mathcal{P}$.

\begin{example}
The operads $Lie$, $Comm$, $Assoc$, the Gerstenhaber and the Poisson
operads, are Koszul. See [GK], Section 4 for a proof.
\end{example}

\begin{definition}
Let $\mathcal{P}$ be a quadratic Koszul operad. A homotopy
$\mathcal{P}$-algebra (or $\mathcal{P}_\infty$-algebra) is an
algebra over the free dg operad $\mathbf{D}(\mathcal{P}^!)$.
\end{definition}
We denote by $\mathcal{P}^*$ the cooperad dual to an operad
$\mathcal{P}$, if all spaces $\mathcal{P}(n)$ are
finite-dimensional. Let $\mathcal{P}$ be a Koszul operad. Then to
define a $\mathcal{P}_\infty$-algebra structure on $X$ is the same
that to define a differential on the free coalgebra
$\mathbb{F}^\vee_{\mathcal{P}^{!*}}(X)$ which is a coderivation of
the coalgebra structure. Any $\mathcal{P}$ algebra is naturally a
$\mathcal{P}_\infty$-algebra.

We denote by $\mathbb{F}_{\mathcal{P}}(V)$ the free algebra over the
operad $\mathcal{P}$, and by $\mathbb{F}^\vee_{\mathcal{P}^*}$ the
free coalgebra over the cooperad $\mathcal{P}^*$. Here we suppose
that all spaces $\mathcal{P}(n)$ are finite-dimensional.

Recall the following statement [GK], Thm. 4.2.5:

\begin{lemma}
Let $\mathcal{P}$ be a Koszul operad, and $V$ a vector space. Let
$X=\mathbb{F}_{\mathcal{P}}(V)$. Then the natural projection
\begin{equation}\label{2.3.1}
(\mathbb{F}^\vee_{\mathcal{P}^{!*}}(X))\to V
\end{equation}
is a quasi-isomorphism.
\end{lemma}

It follows from this statement that any $\mathcal{P}$-algebra $A$
has the following free resolution $\mathcal{R}^\mb(A)$:
\begin{equation}\label{2.3.2}
\mathcal{R}^\mb(A)=(\mathbb{F}_{\mathcal{P}}(\mathbb{F}^\vee_{\mathcal{P}^{!*}}(A),Q_1),Q_2)
\end{equation}
with the natural differentials $Q_1$ and $Q_2$.

Now we define the {\it truncated deformation complex} of the
$\mathcal{P}$-algebra $A$ as $(\mathrm{Der}(\mathcal{R}^\mb(A)),Q)$
where $Q$ comes from the differential in $\mathcal{R}^\mb(A)$. This
deformation complex is naturally a dg Lie algebra with the Lie
bracket of derivations. We have the following statement:

\begin{proposition}
The truncated deformation functor associated with this dg Lie
algebra governs the formal deformations of $A$ as
$\mathcal{P}_\infty$-algebra.
\end{proposition}

\begin{remark}
The word "truncated" means that for the "full" deformation functor
we should take the quotient modulo the inner derivations. Although,
a map $X\to \mathrm{Der}(X)$ is not defined for an arbitrary operad.
Our truncated deformation functor looks like the Hochschild
cohomological complex of $A$ without the degree 0 term $A$.
\end{remark}

The following trick simplifies computations with the deformation
complex.

Any coderivation of the coalgebra
$(\mathbb{F}^\vee_{\mathcal{P}^{!*}}(A),Q_1)$ can be extended to a
derivation of $\mathcal{R}^\mb(A)$ by the Leibniz rule. We have the
following theorem:

\begin{theorem*}
The natural inclusion
\begin{equation}\label{2.3.3}
\mathrm{Coder}(\mathbb{F}^\vee_{\mathcal{P}^{!*}}(A),Q_1)\rightarrow
\mathrm{Der}(\mathcal{R}^\mb(A))
\end{equation}
is a quasi-isomorphism of dg Lie algebras.
\end{theorem*}

It follows from this Theorem and the Proposition above that the dg
Lie algebra
$\mathrm{Coder}(\mathbb{F}^\vee_{\mathcal{P}^{!*}}(A),Q_1)$ governs
the formal deformation of the $\mathcal{P}_\infty$-structure on $A$.
This, of course, can be seen more directly. Indeed, a
$\mathcal{P}_\infty$ structure on $A$ is a differential on
$\mathbb{F}^\vee_{\mathcal{P}^{!*}}(A)$ making latter a dg coalgebra
over $\mathcal{P}^!$. We have some distinguished differential $Q_1$
on it, arisen from the $\mathcal{P}$-algebra structure on $A$. When
we deform it, it is replaced by $Q_\hbar=Q_1+\hbar d_\hbar$ such
that $(Q_1+\hbar d_\hbar)^2=0$. In the first order in $\hbar$ we
have the condition $[Q_1, d_\hbar]=0$, where the zero square
condition is the Maurer-Cartan equation in the corresponding dg Lie
algebra.

\subsection{The main computation in the Tamarkin's theory}\label{section2.4}
Here we compute the deformation cohomology of $T_\poly(V)$, $V$ a
complex vector space, as Gerstenhaber algebra. We prove here Lemma
\ref{section2.2}, and a more general statement.

We start with the following Lemma:
\begin{lemma}
The Gerstenhaber operad $G$ is Koszul. The Koszul dual operad $G^!$
is $G[-2]$.
\begin{proof}
We know that $Lie^!=Comm[-1]$ and $Comm^!=Lie[-1]$. A structure of a
Gerstenhaber algebra on $W$ consists from compatible actions of
$Comm$ and $Lie[1]$ on $W$. The quadratic dual to $Comm$ is
$Lie[-1]$, and the quadratic dual to $Lie[1]$ is $Comm[-2]$.
Therefore, the quadratic dual to $G$ is $G[-2]$. The Koszulity of
$G$ is proven in [GJ], see also [GK] and [H].
\end{proof}
\end{lemma}

Theorem 2.3 gives us a way how to compute the deformation functor
for formal deformations of $T_\poly(V)$ as homotopy Gerstenhaber
algebra. We take the free coalgebra
$\mathbb{F}_{G[-2]^*}^\vee(T_\poly(V))$ over the cooperad $G[-2]^*$
cogenerated by $T_\poly(V)$. It is clear that
\begin{equation}\label{2.4.1}
\mathbb{F}^\vee_{G^*[2]}(T_\poly(V))=S^\mb((\mathbb{F}_{\Lie}
T_\poly(V)[1])[1])[-2]
\end{equation}

The product $\wedge\colon S^2T_\poly(V)\to T_\poly(V)$ and the Lie
bracket $[\ ,\ ]\colon S^2 T_\poly(V)\to T_\poly(V)[-1]$
define two {\it coderivations} of the Gerstenhaber coalgebra
structure on $\mathbb{F}^\vee_{G^*[2]}(T_\poly(V))$; denote them by
$\delta_{\Comm}$ and $\delta_\Lie$, correspondingly.

The deformation complex of $T_\poly(V)$, as of
Gerstenhaber algebra, is then
$\mathrm{Coder}^\mb(\mathbb{F}^\vee_{G^*[2]}(T_\poly(V)))$ endowed
with the differential $d=ad(\delta_{\Comm})+ad(\delta_\Lie)$. We
denote the two summands by $d_\Comm$ and $d_\Lie$, correspondingly.

\begin{theorem*}
The $\mathrm{Aff}(V)$-invariant subcomplex in the ``positive'' deformation complex $(\mathrm{Coder}^+(\mathbb{F}^\vee_{G^*[2]}(T_\poly(V))),
d_\Comm+d_\Lie)$ has all vanishing cohomology.

\begin{remark}
Here the {\it positive} deformation complex means that we exclude the constant coderivations $\Hom(\mathbb{C},T_\poly(V)[2])$. The reason to consider the positive complex is that the constant coderivations do not appear in the formal deformations of $G_\infty$ algebra structure.
\end{remark}

\begin{proof}
First of all we compute the cohomology of all derivations, then concentrate on positive $\mathrm{Aff}$-invariant subcomplex.

Recall in the beginning some tautological facts. The free Gerstenhaber algebra generated by a vector
space $W$ is $S^\mb(\mathbb{F}_\Lie(W[1])[-1])$. The cofree Gerstenhaber coalgebra cogenerated by $W$ is
$S^\mb(\mathbb{F}^\vee_\Lie(W[-1])[1])$. Finally, the cofree $G^*[2]$-coalgebra cogenerated by $W$ is
$S^\mb(\mathbb{F}^\vee_\Lie(W[1])[1])[-2]$.

We deal with the coderivations of the  cofree coalgebra $S^\mb(\mathbb{F}^\vee_\Lie(W[1])[1])[-2]$, and
they are defined uniquely by their restrictions to cogenerators
which may be arbitrary. Therefore, we need to compute the cohomology
of the complex
\begin{equation}\label{2.4.2}
\Hom_\mathbb{C}(S^\mb((\mathbb{F}_\Lie (T_\poly(V)[1]))[1]),
T_\poly(V))[2]
\end{equation}
This is a bicomplex with the differentials $d_\Comm$ and $d_\Lie$.
We use the spectral sequence of the bicomplex which computes firstly the
cohomology of $d_\Comm$. We leave to the reader the simple check that this spectral sequence converges to the (associated graded of) cohomology of the total complex.

Compute the first term of the spectral sequence. When we take in the account
the only differential $d_\Comm$, the deformation complex
$\Hom_\mathbb{C}(S^\mb((\Free_\Lie (T_\poly(V)[1]))[1]),
T_\poly(V))[2]$ is a direct sum of {\it complexes}:

\begin{equation}\label{2.4.3}
\begin{aligned}
\ &\Hom_\mathbb{C}(S^\mb((\Free_\Lie (T_\poly(V)[1]))[1]),
T_\poly(V))[2]=\\
&T_\poly(V)[2]\bigoplus_{k\ge 1}(\Hom( S^k((\Free_\Lie T_\poly(V)[1])[1]),
T_\poly(V)),d_\Comm)[2]
\end{aligned}
\end{equation}

\begin{lemma}
Let $k\ge 1$. The cohomology of the complex $(\Hom( S^k((\Free_\Lie T_\poly(V)[1])[1]),
T_\poly(V)),d_\Comm)[2]$ is $S^k_\mathcal{O}(\mathrm{Vect}(T^*[-1]V))[-2k+2]$. Here $S^k_\mathcal{O}(\mathrm{Vect}(T^*[-1]V))$ is (the sections of) the $k$-th symmetric power of the vector bundle of vector fields on the space $T^*[-1]V$. (Recall that $T_\poly(V)$ is the functions on $T^*[-1]V$).
\end{lemma}
\proof{}
We only ``explain'' the statement, the complete proof will appear somewhere.
The complex $(\Hom( S^k((\Free_\Lie T_\poly(V)[1])[1]),
T_\poly(V)),d_\Comm)[2]$ ``starts'' with the term $\Hom(S^k (T_\poly(V)), T_\poly(V))[-2k+2]$ (we take only the generators of the free Lie (co)algebra). The differential $d_\Comm$ in this term is
\begin{equation}
\begin{aligned}
\ &(d_\Comm\Psi)(\gamma_1\cdot\dots\cdot\gamma_{k+1})=\\
&\mathrm{Symm}\Psi((\gamma_1\wedge\gamma_2)\cdot\gamma_3\cdot\dots\cdot\gamma_{k+1})\mp \\
&\mathrm{Symm}\bigl(\gamma_1\wedge\Psi(\gamma_2\cdot\dots\cdot\gamma_{k+1})\pm \gamma_2\wedge\Psi(\gamma_1\cdot\gamma_3\cdot\dots\cdot\gamma_{k+1})\bigr)
\end{aligned}
\end{equation}
The kernel of this differential is exactly the answer given in the statement of the Lemma (this is clear). The more non-trivial is to show that the ``higher'' cohomology vanishes.
\endproof
Thus, in the term $E_1$ we have $\bigoplus_{k\ge 0}S^k_\mathcal{O}(\mathrm{Vect}(T^*[-1]V))[-2k+2]$.

Now consider the differential $d_\Lie$ acting on $E_1$. The Schouten bracket is an element of $S^2_\mathcal{O}\mathrm{Vect}(T^*[-1]V)[-1]$.
One easily sees that the cohomology belongs to $T_\poly(V)[2]\subset \bigoplus_{k\ge 0}S^k_\mathcal{O}(\mathrm{Vect}(T^*[-1]V))[-2k+2]$ (the summand for $k=0$). This cohomology is 1-dimensional and is represented by a constant function in $T^0_\poly(V)[2]$.
We conclude that the cohomology of the full deformation complex \eqref{2.4.2} is 1-dimensional and is concentrated in degree -2.

Now consider the $\Aff$-invariant subcomplex of the full deformation complex.
More precisely, we compute the cohomology of
\begin{equation}\label{ref1}
\Hom_\mathbb{C}(S^\mb((\Free_\Lie (T_\poly(V)[1]))[1]),
T_\poly(V))^\Aff[2]
\end{equation}
From a very general point of view, the terms of our complex are $\Aff$-invariants in $\Hom(V^{\otimes a}\otimes V^{*\otimes b},V^{\otimes c}\otimes V^{*\otimes d})$. The Lie group $\mathrm{GL}(n)$ acts in the natural way, and the shift $x\mapsto x+(a_1,\dots,a_n)$ acts trivially on all $V$ factors, and by shifts $x_i\mapsto x_i+a_i$ on all $V^*$ factors.

We know all $\mathrm{GL}(n)$-invariants from the main theorem of invariant theory. They are constructed from the following 4 elementary operations. These 4 operations are: $id: V\to V$, $id: V^*\to V^*$, $V\otimes V^*\to\mathbb{C}$, and $\mathbb{C}\to V\otimes V^*$. From these 4 operations, only the last one, $\mathbb{C}\to V\otimes V^*$, is not $\Aff$-invariant. On the other hang, the group $\mathrm{GL}(n)$ acts reductively, and the cohomology of the complex is equal to the cohomology of its $\mathrm{GL}(n)$-invariant part. 

On the other hand, due to the symmetrization conditions, the invariant $c\colon\mathbb{C}\to V\otimes V^*$ can be applied only 0 or 1 times. 
So schematically as a vector space the space of $\mathrm{GL}(n)$-invariants of \eqref{2.4.3} is $K^\mb\oplus c\cdot K^\mb$ where $K^\mb$ is the space of $\Aff$-invariants of \eqref{2.4.3}. One easily sees that this decomposition agrees with the complex structure. We conclude that the cohomology of \eqref{ref1} is equal to the $\Aff$-invariants of the cohomology of
\begin{equation}\label{ref2}
\Hom_\mathbb{C}(S^\mb((\Free_\Lie (T_\poly(V)[1]))[1]),
T_\poly(V))^{\mathrm{GL}(n)}[2]
\end{equation}
As was already mentioned above, the latter cohomology is equal to the $\Aff$-invariants of the cohomology of \eqref{ref1}, because the group $\mathrm{GL}(n)$ acts reductively. 

We conclude that the cohomology of \eqref{ref2} is equal to $\mathbb{C}[2]$.

The last step is to compute the cohomology of the positive subcomplex. This is easy to do. In the term $E_2$ one has $(T_\poly(V)/\mathbb{C})[0]$, and after taking of the $\Aff$-invariants, we get 0.

Theorem is proven.
\end{proof}
\end{theorem*}

\subsection{The final point: relation with the Etingof-Kazhdan quantization}\label{section2.5}
The remaining part of the Tamarkin's proof of Kontsevich formality
goes as follows.

One firstly proves the {\it Deligne conjecture} that there is a
homotopy Gerstenhaber algebra structure on the Hochschild
cohomological complex $\Hoch^\mb(A)$ of any associative algebra such
that it induces the Schouten bracket and the wedge product on the
cohomology. This is the only transcendental step of the
construction, this structure, as it is defined in [T], depends on a
choice of Drinfeld associator [D].

We apply this fact for $A=S(V^*)$, $V$ a vector space over
$\mathbb{C}$. One can push-forward (given by the "Massey
operations") this $G_\infty$ structure from $\Hoch^\mb(S(V^*))$ to
its cohomology $T_\poly(V)$. Then we get two $G_\infty$ structures
on $T_\poly(V)$, the first is given from the Schouten bracket and
the wedge product of polyvector fields, the second is the above
pushforward. Moreover, one can introduce a formal parameter $\hbar$
to the pushforward, such that the original one is given when
$\hbar=1$. Then for $\hbar=0$ we get the Schouten structure: it
follows from the compatibility of the Deligne's conjecture
$G_\infty$ structure with the one on the cohomology. Thus we get a
formal deformation of the classical Gerstenhaber algebra structure
on $T_\poly(V)$. This deformation is clearly
$\mathrm{Aff}(V)$-invariant. By Theorem 2.4, infinitesimally all
such deformations are trivial; therefore, they are trivial globally.
This concludes the Tamarkin's proof.

In this Subsection we explain the Tamarkin's proof of the Deligne
conjecture, based on the Etingof-Kazhdan quantization.

Recall the definitions of the dg operads $B_\infty$ and $B_\Lie$. By
definition, a vector space $X$ is an algebra over the operad
$B_\infty$ if there is a structure of a dg associative bialgebra on
the cofree coalgebra $\mathbb{F}^\vee_\Assoc(X[1])$ such that the
coalgebra structure coincides with the given one. This definition
leads to the following data (see [H], Section 5 and [GJ], Section 5
for more detail):
\begin{itemize}
\item[(1)] a differential $d\colon X[1]^{\otimes n}\to X[1]^{\otimes
m}$ of degree 1, $m,n\ge 1$ being a differential of the free
coalgebra structure is uniquely defined by the projections to the
cogenerators. We denote them $m_n\colon X[1]^{\otimes n}\to X[2]$;
\item[(2)] the algebra structure, it is also given by the projection
to $X[1]$. These are maps $m_{pq}\colon X[1]^{\otimes p}\otimes
X[1]^{\otimes q}\to X[1]$, or $m_{pq}\colon X^{\otimes p}\otimes
X^{\otimes q}\to X[1-p-q]$.
\end{itemize}
These data should define three series of equations: the first come
from the associativity of $\{m_{pq}\}$, the second come from the
fact that $d$ is a derivation of the algebra structure, and the
third comes from the condition $d^2=0$. These equations define a
very complicated operad $B_\infty$.

It is a remarkable and surprising result of Getzler-Jones [GJ],
Section 5, that $X=\Hoch^\mb(A)$, $A$ an arbitrary associative
algebra, is an algebra over the operad $B_\infty$. This structure is
defined as follows:
\begin{itemize}
\item[(1)] $m_1$ is the Hochschild differential;
\item[(2)] $m_2$ is the cup-product on $\Hoch^\mb(A)$;
\item[(3)] $m_i=0$ for $i\ge 3$;
\item[(4)] $m_{1k}(f\otimes g_1\otimes \dots\otimes g_k)$ is the
brace operation $f\{g_1,\dots,g_k\}$ defined below;
\item[(5)] $m_{ak}=0$ for $a\ge 2$.
\end{itemize}
Now is the definition of the braces due to Getzler-Jones. It is
better to describe it graphically, as is shown in Figure 2.
\sevafigc{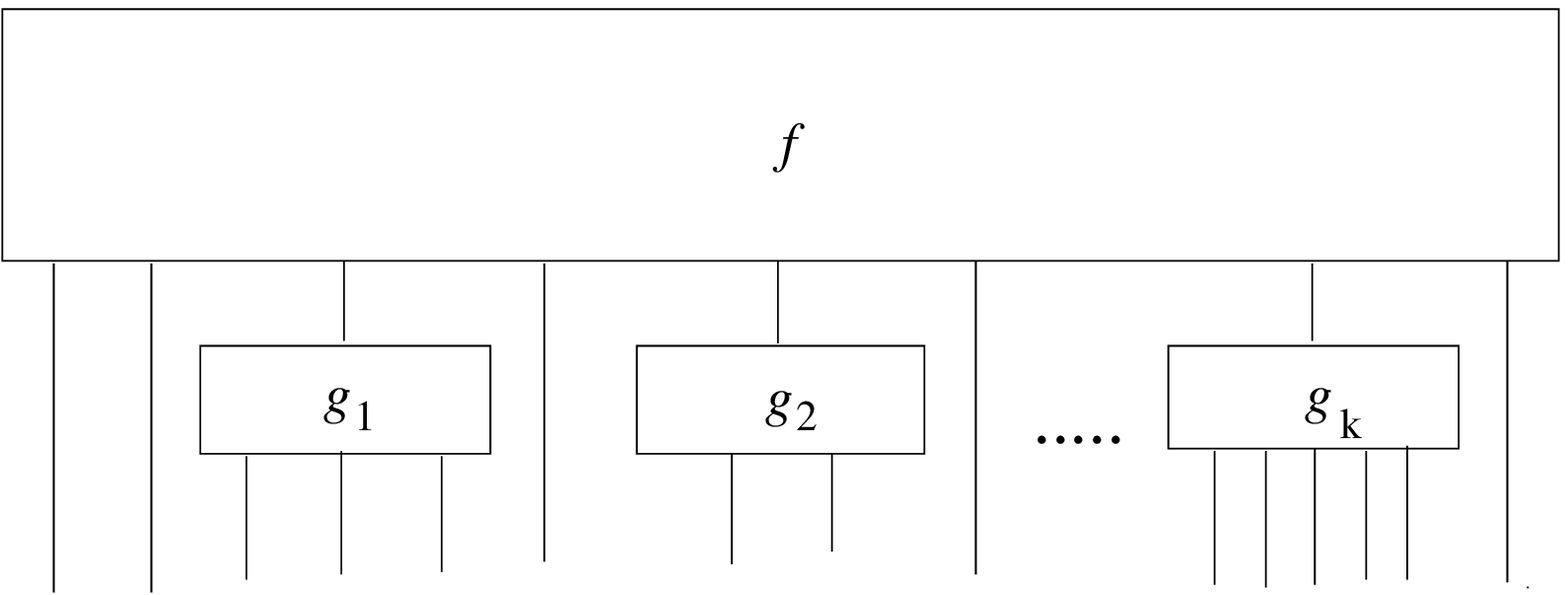}{80mm}{0}{The brace operation: we insert
$g_1,\dots, g_k$ into arguments of $f$, preserving the order of
$g_1,\dots, g_k$, with the natural sign, and take the sum over all
possible insertions}

Let us emphasize again that it is a highly non-evident fact, proven
by a direct computation, that in this way we make a $B_\infty$
algebra structure on $\Hoch^\mb(A)$.

The role of this construction is that the cohomology operad of the
dg operad $B_\infty$ is equal to the Gerstenhaber operad $G$
(probably, even to prove this fact we need the Etingof-Kazhdan
quantization). So the idea is to prove that there is a
quasi-isomorphism of operads $G\rightsquigarrow B_\infty$, and then
we can consider $\Hoch^\mb(A)$ as $G_\infty$ algebra for any
associative algebra $A$. The only trancendental step in the
Tamarkin's construction is this quasi-isomorphism of operads
$G\rightsquigarrow B_\infty$.

Technically it is done as follows. Introduce some operad $B_\Lie$ as
follows. A vector space $X$ is an algebra over the operad $B_\Lie$
if there is a dg Lie bialgebra structure on the free Lie coalgebra
$\mathbb{F}^\vee_\Lie(X[1])$ generated by $X[1]$ such that the Lie
coalgebra structure coincides with the given one.

The operad $B_\Lie$ is also quasi-isomorphic to the Gerstenhaber
operad $G$. Moreover, there is a simple construction of a
quasi-isomorphism $G_\infty\to B_\Lie$, as follows.

Let $Y$ be an algebra over $B_\Lie$. This means that there is a Lie
dg bialgebra structure on the free coalgebra
$\g=\mathbb{F}^\vee_\Lie(Y[1])$. In particular, $\g$ is a Lie
algebra, and this defines a differential on the Lie chain complex
$\mathbb{F}^\vee_{\mathrm{Com}}(\g[1])$. Thus we get a differential
on the free Gerstenhaber coalgebra $\mathbb{F}^\vee_{G^\vee}(Y[1])$
cogenerated by $Y[1]$, which by definition means that $Y$ is a
$G_\infty$-algebra. This assignment is functorial, and therefore
gives a map of operads $G_\infty\to B_\Lie$, which easily checked to
be a quasi-isomorphism.

So, the conclusion is that the operad $B_\Lie$ can be connected to
the Gerstenhaber operad in a very simple way, and now we should
connect the operad $B_\infty$ with the operad $B_\Lie$.

This is given exactly by the Etingof-Kazhdan (de)quantization (see
[T1] and [H], Section 7 for detail). The Etingof-Kazhdan
dequantization is applied in a sense to
$\mathbb{F}^\vee_\Assoc(\Hoch^\mb(A)[1])$ which is an associative
bialgebra by the Getzler-Jones braces' construction.

\begin{remark}
Let $P$ be a Poisson algebra. Its Poisson complex $\Pois^\mb(P)$ is
defined as the dg space of coderivations of the free coalgebra over
the dual cooperad $\mathcal{P}^{!*}[1]$ by the space $P[1]$. This
space of coderivations is naturally equipped with a differential
$d_\Pois$ arising from the Poisson bracket and the product on $P$.
The author thinks that for any Poisson algebra $P$ the Poisson
complex $\Pois^\mb(P)$ is an algebra over the operad $B_\Lie$. This
structure is defined exactly by some generalization of the braces
construction. Now, if $P=S(V^*)$ be the Poisson algebra with zero
bracket, then by Getzler-Jones $\Hoch^\mb(P)$ is a $B_\infty$
algebra, and $\Pois^\mb(P)$ is a $B_\Lie$ algebra. The author thinks
that some the most natural Etingof-Kazhdan dequantization gives from
the associative bialgebra $\mathbb{F}^\vee_\Assoc(\Hoch^\mb(P)[1])$
the Lie bialgebra $\mathbb{F}^\vee_\Lie(\Pois^\mb(P)[1])$. So far,
the author does not know any direct proof of the last fact.
\end{remark}

\subsection*{Appendix}
Here we explain a construction of the Kontsevich formality morphism
over $\mathbb{Q}$. The usual construction uses the Drinfeld's
associator over $\mathbb{Q}$ and the Tamarkin's theory. This
associator is not given by an explicit formula, it is constructed in
[D] by proving that all associators over $\overline{\mathbb{Q}}$
form a torsor over the Grothendieck-Teichm\"{u}ller group. The
Knizhnik-Zamolodchikov associator gives an example of associator
over $\mathbb{C}$; therefore, there exists an associator over
$\overline{\mathbb{Q}}$. It proves that this torsor is trivial over
$\overline{\mathbb{Q}}$. Then the torsor is trivial also over
$\mathbb{Q}$, because the Grothendieck-Teichm\"{u}ller group is
pro-unipotent.

Here we propose a different proof, which seems to be more
constructive. This approach seems to be new.

It follows from the previous results that if we succeed to construct
a quasi-isomorphism of operads $B_\infty\to G$ over $\mathbb{Q}$, we
will be done.

Consider the operad $B_\infty$. It is a dg operad. All dg operads
form a closed model category because they are algebras over some
universal colored operad. In particular, there is a homotopy operad
structure on the cohomology of $B_\infty$, given by a kind of
"Massey operations". This homotopy operad structure clearly is
defined over $\mathbb{Q}$. To construct it explicitly, we should
firstly split $B_\infty$ as a complex into a direct sum of its
cohomology and a contractible complex and, secondly, to contract
this complex explicitly by a homotopy. It is clear that these two
steps can be performed over $\mathbb{Q}$. Now we have two homotopy
operad structures on $G$: one is the Gerstenhaber operad, and
another one is given by the Massey operations. Moreover, there is a
formal family of homotopy dg operads depending on $\hbar$ whose
value at $\hbar=1$ is the "Massey" homotopy operad structure on $G$,
and whose value at $\hbar=0$ is the Gerstenhaber operad.

We know from the Tamarkin's theory described above that this
deformation is trivial over $\mathbb{C}$, because the "Massey"
homotopy operad is quasi-isomorphic to $B_\infty$ by construction,
which is quasi-isomorphic over $\mathbb{C}$ to the Gerstenhaber
operad by Section 2.5. We are going to prove that this formal
deformation is trivial also over $\mathbb{Q}$.

For this it is sufficient to prove that infinitesimally this
deformation is trivial over $\mathbb{Q}$ at each $0\le \hbar\le 1$.
Consider a resolution $\mathcal{R}^\mb(G)$ of the Gerstenhaber
operad over $\mathbb{Q}$; as $G$ is a Koszul operad, we can take its
Koszul resolution. Consider the truncated deformation complex
$\mathbb{D}^+(G)=\mathrm{Der}(\mathcal{R}^\mb(G))$. We need to prove
that all infinitesimal deformations give trivial classes in
$H^1(\mathbb{D}^+(G))$.

It is probably not true that $H^1(\mathbb{D}^+(G))=0$, it would be
very unexpected, because Tamarkin imbeded in [T3] the
Grothendieck-Teichm\"{u}ller Lie algebra into
$H^0(\mathbb{D}^+(G))$. But we know that all classes are trivial in
$H^1(\mathbb{D}^+(G),\mathbb{C})$ from the Tamarkin's theory. As the
complex $\mathbb{D}^+(G)$ is defined over $\mathbb{Q}$, we have that
the natural map
\begin{equation}\label{2.6.1}
H^\mb(\mathbb{D}^+(G),\mathbb{Q})\hookrightarrow
H^\mb(\mathbb{D}^+(G)\otimes_\mathbb{Q}\mathbb{C})=H^\mb(\mathbb{D}^+(G),\mathbb{C})
\end{equation}
is an embedding.

Therefore, all our infinitesimal classes are trivial over
$\mathbb{Q}$, and we get that the global formal deformation is
trivial over $\mathbb{Q}$.

We think that this speculation is as explicit as it can give some
explicit formulas for the Kontsevich's formality over $\mathbb{Q}$.
We are going to describe it in detail in the sequel.

\section{Two infinite-dimensional varieties (and a morphism between them)}
In this Section we associate with each quasi-isomorphism of operads
$\Theta\colon G_\infty\to B_\infty$ defined up to homotopy an
$L_\infty$ morphism $\U(\Theta)\colon T_\poly(V)\to \Hoch^\mb(V)$
(for any vector space $V$) defined up to homotopy. We show that this
$L_\infty$ morphism is given by a universal formula, that is by
prediction of some weights to all Kontsevich graphs from [K97], but
our graphs may contain simple loops. The image of the map
$\Theta\mapsto \U(\Theta)$ gives that $L_\infty$ morphisms for
deformation quantizations associated with which we can prove our
Koszul duality Theorem.

\subsection{Few words about homotopy}
Starting from now, we often use the word "homotopy" in the context
like "homotopical map of dg operads" or "homotopical $L_\infty$
morphisms". Here are some generalities on this.

The Quillen's formalism of closed model categories [Q] gives a tool
for the inverting of quasi-isomorphisms in a non-abelian case. Let
$\mathcal{O}$ be an operad. Consider the category $DGA(\mathcal{O})$
of dg algebras over $\mathcal{O}$. We want to construct a universal
category in which the quasi-isomorphisms in $DGA(\mathcal{O})$ are
invertible. This category can be constructed for any operad
$\mathcal{O}$ and is called the homotopy category of
$DGA(\mathcal{O})$, because the category $DGA(\mathcal{O})$ admits a
closed model structure in which the weak equivalences coincide with
the quasi-isomorphisms [H2]. There are several constructions of this
category, but all them are equivalent due to the universal property
with respect to the localization by quasi-isomorphisms. In Section 5
we recall a very explicit construction in the case when
$\mathcal{O}$ is a Koszul operad.

Contrary, the dg algebras over a PROP do not form a closed model
category (the Quillen's Axiom 0 that the category admits all finite
limits and colimits fails in this case; for example we do not know
what is a free algebra over a PROP). Therefore, for dg algebras over
a PROP any construction of the homotopical category (to the best of
our knowledge) is not known.

On the other hand, all dg operads form a closed model category as
algebras over the universal colored operad, therefore, for dg
operads the Quillen's construction works.

In the sequel we will skip some details concerning that homotopical
maps of operads induce homotopical maps of dg algebras, avoiding to
enlarge this already rather long paper.

Only what we need to know is that the homotopical category is
unique, and in the final step we use a particular construction of it
for dg Lie algebras in Section 5, appropriate for our needs.

\subsection{The Kontsevich's variety $\mathfrak{K}$}
The Kontsevich's variety $\mathfrak{K}$ is the variety of all {\it
universal} $L_\infty$ quasi-isomorphisms
$T_\poly(V)\to\Hoch^\mb(S(V^*))$ defined for all vector spaces $V$.
Any such universal $L_\infty$ morphism is by definition given by
prediction of some complex weights $W_\Gamma$ to all Kontsevich
graphs $\Gamma$ from [K97] possibly with simple loops and not
connected. These $W_\Gamma$ are subject to some quadratic equations,
arising from the $L_\infty$ condition. The first
Hochschild-Kostant-Rosenberg graph has the fixed weight, as in the
Kontsevich's paper [K97]. This variety is not empty as is proven in
[K97]. The homotopies acts by gauge action (see Section 5).

\subsection{The Tamarkin's variety $\mathfrak{T}$}
The Tamarkin's variety in our strict sense is formed from all
quasi-isomorphisms of operads $G_\infty\to B_\infty$ which are
identity on cohomology, modulo homotopies. As $G_\infty$ is a free
dg operad, any such map is uniquely defined by the generators
$G[-2]$. So, a map of operads $G_\infty\to B_\infty$ is given by a
map of vector spaces $G[-2]\to B_\infty$ subject to some quadratic
relations arose from the compatibility with the differentials in the
dg operads. This variety is not empty because we have constructed in
Section 2, following [T1], such a particular quasi-isomorphism.

In a wider setting, one can consider $Op_\infty$ maps of dg operads
$G_\infty\to B_\infty$, but we do not do this.

\subsection{A map $\mathfrak{X}\colon\mathfrak{T}\to \mathfrak{K}$}
Suppose a point $t$ of the Tamarkin's variety $\mathfrak{K}$ is
fixed. Then the Hochschild complex $\Hoch^\mb(A)$ of any algebra $A$
has a fixed structure of homotopy Gerstenhaber algebra (fixed modulo
homotopy). Consider the case $A=S(V^*)$ for some vector space $V$.
Then we get, as is explained in Section 2, two structures of
$G_\infty$ algebra on $T_\poly(V)$ which are specifications of some
formal deformation at $\hbar=0$ and $\hbar=1$. Then they should
coincide, up to a homotopy, because the first cohomology
$H^1(\mathrm{Coder}(\mathbb{F}^\vee_{G^*[2]}(T_\poly(V))))=T^2_\poly(V)$
by Theorem 2.4, and there is no $\mathrm{Aff}(V)$-invariant classes
(but our formal deformation is $\mathrm{Aff}$-invariant).

Thus we get a map of $G_\infty$ algebras $\mathfrak{X}^0(t)\colon
T_\poly(V)\to\Hoch^\mb(S(V^*))$, where $T_\poly(V)$ is considered
with the standard Schouten-Nijenhuis Gerstenhaber structure, and the
$G_\infty$ structure on $\Hoch^\mb(S(V^*))$ depends on the point
$t\in \mathfrak{K}$. Then we restrict it to the Lie operad and get
an $L_\infty$ map $\mathfrak{X}(t)\colon
T_\poly(V)\to\Hoch^\mb(S(V^*))$ which is an $L_\infty$ morphism for
the standard Lie structures on $T_\poly(V)$ and $\Hoch^\mb(S(V^*))$,
and this $L_\infty$ morphism is defined up to homotopy. This is the
construction of the map $\mathfrak{X}$. It is, although, not proven
yet that $\mathfrak{X}(t)$ is defined uniquely up to a homotopy.

\begin{lemma} For a fixed $V$, the
$L_\infty$ morphism $\mathfrak{X}(t)$ is uniquely defined up to a
homotopy.
\begin{proof}
Suppose there are two different $G_\infty$ morphisms, for the same
fixed $G_\infty$ structure on $\Hoch^\mb(S(V^*))$. Then we can get
defined up to a homotopy $G_\infty$ quasi-automorphism of
$T_\poly(V)$. It has the identity first Taylor component by the
constructions (a point $t\in \mathfrak{T}$ is defined such we get
the canonical Gerstenhaber structure on the cohomology of
$\Hoch^\mb(A)$ for any $A$). Therefore, the logarithm of this
automorphism is well defined and gives a $G_\infty$ derivation of
$T_\poly(V)$. Now we use the computation of Theorem 2.4 for 0-th
cohomology:
$H^0(\mathrm{Coder}(\mathbb{F}^\vee_{G^*[2]}(T_\poly(V))))=T^1_\poly(V)$
is the vector fields. Again, there are no $\mathrm{Aff}(V)$-ivariant
vector fields. Therefore, our derivation is inner. But then it is
zero, because any inner derivation acts non-trivially on the first
Taylor component which is fixed to be identity. Thus, we have proved
that $\mathfrak{X}(t)$ is well-defined up to homotopy as $G_\infty$
map, and therefore the same is true for its $L_\infty$ part. (Compare with Lemma in the end of Section 4.4).
\end{proof}
\end{lemma}

Now we prove the following almost evident corollary of the previous
Lemma:

\begin{theorem*}
The $L_\infty$ morphism $\mathfrak{X}(t)$ is a universal $L_\infty$
morphism, that is, it is given by prediction of some weights to all
Kontsevich graphs, possibly non-connected and with simple loops, and
these weights up to a homotopy do not depend on the vector space
$V$.
\begin{proof}
Let $W\subset V$ be a subspace. Then we can decompose $V=W\oplus
W^\bot$, and a $G_\infty$ structure on $\Hoch^\mb(S(V^*))$ defines a
$G_\infty$ structure on $\Hoch^\mb(S(W^*))$. Clearly (because the
$G_\infty$ structures are $\gl(V)$-invariant) it is, up to a
homotopy, the structure on $\Hoch^\mb(S(W^*))$ one gets from the map
of operads $t\colon G_\infty\to B_\infty$. We have then two
definitions of $L_\infty$ morphisms
$T_\poly(W)\to\Hoch^\mb(S(W^*))$: one is the direct
$\mathfrak{X}(t)_W$, and the second one is the restriction to $W$ of
$\mathfrak{X}(t)_V$. They coincide up to a homotopy by the Lemma
above, because the two $G_\infty$ structures on $\Hoch^\mb(S(W^*))$
are the same up to a homotopy. The remaining part of the Theorem
(that the $L_\infty$ morphism is given by a universal formula though
Kontsevich graphs) follows from the $\gl(V)$ invariance of it.
\end{proof}
\end{theorem*}

Is is not known if the map $t\mapsto \mathfrak{X}(t)$ is surjective,
even when we allow $t$ to be an $Op_\infty$ map of dg operads
$G_\infty\to B_\infty$. Our Main Theorem of this paper about the
Koszul duality holds for the star-product obtained from any
$L_\infty$ morphism in the image of $\mathfrak{X}$,
$\U=\mathfrak{X}(t)$, by the usual formula
\begin{equation}\label{eq3.4.1}
f\star g=f\cdot g+\hbar\U_1(\alpha)(f\otimes
g)+\frac12\hbar^2\U_2(\alpha,\alpha)(f\otimes g)+\dots
\end{equation}
where $\alpha$ is a (quadratic) Poisson bivector field on $V$.

\section{Koszul duality and dg categories}

\subsection{Some generalities on dg categories}
We give here some basic definitions on dg categories. We define only
the things we will directly use, see [Kel3] for much more detailed
and sophisticated overview.

A {\it differential graded (dg) category} $\mathcal{A}$ over a field
$k$ is a category, in which the sets of morphisms $\mathcal{A}(X,Y)$
between any two objects $X,Y\in \Ob(\mathcal{A})$ are $k$-linear dg
spaces (complexes of $k$-vector spaces) such that the compositions
are defined as maps $\mathcal{A}(Y,Z)\otimes
\mathcal{A}(X,Y)\to\mathcal{A}(X,Z)$ for any
$X,Y,Z\in\Ob(\mathcal{A})$ which are {\it maps of complexes}. In the
last condition we regard $\mathcal{A}(Y,Z)\otimes \mathcal{A}(X,Y)$
as a complex with the differential defined by the Leibniz rule
\begin{equation}\label{eqdgcat1}
d(f\otimes g)=(df)\otimes g+(-1)^{\deg f}f\otimes (dg)
\end{equation}
It is clear that a differential graded category with one object is
just a differential graded associative algebra. Then dg categories
can be considered as "dg algebras with many objects".

For dg algebras we have a definition when a map $F\colon A^\mb\to
B^\mb$ is a {\it quasi-isomorphism}: it means that the map $F$ is a
map of algebras and induces an isomorphism on cohomology. Such a map
in general is not invertible, it can be inverted only as an
$A_\infty$ map.

What should be a definition of {\it a quasi-isomorphism for dg
categories}?

We say that a functor $F\colon\mathcal{A}\to\mathcal{B}$ between two
dg categories $\mathcal{A}$ and $\mathcal{B}$ is a {\it
quasi-equivalence} if it is a functor, which is $k$-linear on
morphisms (and, as such, preserves tensor compositions of morphisms), induces an equivalence on the level of cohomology, and
is essentially surjective. The last
condition means that for a dg category $\mathcal{A}$ we can consider
the category $H^\mb(\mathcal{A})$ with the same objects, and with
$(H^\mb\mathcal{A})(X,Y)=H^\mb(\mathcal{A}(X,Y))$. Then a
quasi-equivalence is not invertible in general, but it can be
inverted as an $A_\infty$ quasi-equivalence between two dg
categories. We will not use this concept directly, and we refer to
the reader to give the definition.

Now if we have a dg algebra, we know what is the cohomological
Hochschild complex of it. It governs the $A_\infty$ deformations of
the dg algebra. It is possible to define the {\it Hochschild
cohomological complex of a dg category}. This will be in a sense the
main object of our study in this paper for some particular dg
category, namely, for the B.Keller's dg category introduced in the
next Subsection. Let us give the definition of it.

At first, it is the total {\it product} complex of a bicomplex. The
vertical differential will be the inner differential appeared from
the differential on $\mathcal{A}(X,Y)$ for any pair
$X,Y\in\Ob\mathcal{A}$. The horizontal differential will an analog
of the Hochschild cohomological differential. The columns have
degrees $\ge 0$. In degree 0 we have
\begin{equation}\label{eqdgcat5}
\Hoch^{*0}(\mathcal{A})=\prod_{X\in\Ob\mathcal{A}}\mathcal{A}(X,X)
\end{equation}
and in degree $p\ge 1$
\begin{equation}\label{eqdgcat6}
\Hoch^{*p}(\mathcal{A})=\prod_{X_0,X_1,\dots,X_p\in\Ob\mathcal{A}}\Hom_k(\mathcal{A}(X_{p-1},
X_p)\otimes \mathcal{A}(X_{p-2},
X_{p-1})\otimes\dots\otimes\mathcal{A}(X_0,X_1),\mathcal{A}(X_0,X_p))
\end{equation}
where the product is taken over all chains of objects
$X_0,X_1,\dots,X_p\in\Ob\mathcal{A}$ of length $p+1$.

The Hochschild differential $d_\Hoch\colon
\Hoch^{*p}(\mathcal{A})\to\Hoch^{*,p+1}(\mathcal{A})$ is defined in
the natural way. Let us note that even if a cochain
$\Psi\in\Hoch^{*p}(\mathcal{A})$ is non-zero only for a single chain
of objects $X_0,X_1,\dots,X_p$, its Hochschild differential
$d_\Hoch\Psi$ in general is non-zero on many other chains of
objects. Namely, at first it may be nonzero for on any chain
\begin{equation}\label{eqdgcat7}
X_0,\dots, X_i,Y, X_{i+1},\dots, X_p \ \text{for}\ 0\le i\le p-1
\end{equation}
such that there are nonzero compositions
$\mathcal{A}(Y,X_{i+1})\otimes\mathcal{A}(X_i,Y)\to\mathcal{A}(X_i,X_{i+1})$
(this is corresponded to the regular terms in the Hochschild
differential), and on the chains
\begin{equation}\label{eqdgcat8}
Z_-,X_0,\dots,X_p
\end{equation}
and
\begin{equation}\label{eqdgcat9}
X_0,\dots, X_p,Z_+
\end{equation}
such that there are non-zero compositions
$\mathcal{A}(X_0,X_p)\otimes\mathcal{A}(Z_-,X_0)\to\mathcal{A}(Z_-,X_p)$
and
$\mathcal{A}(X_p,Z_+)\otimes\mathcal{A}(X_0,X_p)\to\mathcal{A}(X_0,Z_+)$
(this is corresponded to the left and to the right boundary terms in
the Hochschild differential, correspondingly).

The Hochschild cohomological complex of a dg category $\mathcal{A}$
is a dg Lie algebra with the direct generalization of the
Gerstenhaber bracket. The solutions of the Maurer-Cartan equation of
$\Hoch^\mb(\mathcal{A})\otimes k[[\hbar]]$ give the formal
deformations of the dg category $A$ as $A_\infty$ category.

\subsection{The B.Keller's dg category $\cat(A,B,K)$}
Introduce here the main object of our story---the B.Keller's dg
category $\cat(A,B,K)$. Here $A$ and $B$ are two dg associative
algebras, and $K$ is $B$-$A$-bimodule. The dg category
$\mathcal{A}=\cat(A,B,K)$ has two objects, called say $a$ and $b$,
such that $\mathcal{A}(a,a)=A$, $\mathcal{A}(b,b)=B$,
$\mathcal{A}(a,b)=0$, $\mathcal{A}(b,a)=K$. Only what we need from
$K$ to define such a dg category is a structure on $K$ of
differential graded $B$-$A$-bimodule.

Consider in details the Hochschild complex of the category
$\cat(A,B,K)$. It contains as subspaces $\Hoch^\mb(A)$ and
$\Hoch^\mb(B)$, the usual Hochschild cohomological complexes of the
algebras $A$ and $B$, and also it contains the subspace
\begin{equation}\label{eqdgcat10}
\Hoch^\mb(B,K,A)=\sum_{m_1,m_2\ge 0}\Hom(B^{\otimes m_1}\otimes
K\otimes A^{\otimes m_2}, K)[-m_1-m_2]
\end{equation}
As a graded space,
\begin{equation}\label{eqdgcat11}
\Hoch^\mb(\cat(A,B,K))=\Hoch^\mb(A)\oplus\Hoch^\mb(B)\oplus\Hoch^\mb(B,K,A)
\end{equation}
but certainly it is {\it not} a direct sum of {\it subcomplexes}.
Namely, $\Hoch^\mb(B,K,A)$ is a subcomplex of
$\Hoch^\mb(\cat(A,B,K))$, but $\Hoch^\mb(A)$ and $\Hoch^\mb(B)$ are
{\it not}. There are well-defined projections
$p_A\colon\Hoch^\mb(\cat(A,B,K))\to\Hoch^\mb(A)$ and
$\Hoch^\mb(\cat(A,B,K))\to\Hoch^\mb(B)$.

The Hochschild component of the total differential acts like this:

\begin{equation}\label{eqdgcat12}
\xymatrix{& X_3 \ar@(ur,ul)[]_{d_\Hoch^K}\\
 X_1 \ar@(ul,dl)[]_{d_\Hoch^A} \ar[ur]^{d_\Hoch^{AK}} &&X_2 \ar@(ur,dr)[]^{d_\Hoch^B} \ar[ul]_{d_\Hoch^{BK}} }
\end{equation}

\bigskip
where $X_1=\Hoch^\mb(A)$, $X_2=\Hoch^\mb(B)$,
$X_3=\Hoch^\mb(B,K,A)$.

In [Kel1], Bernhard Keller poses the following question: what is a
sufficient condition on the triple $(A,B,K)$ which would guarantee
that the projections
$p_A\colon\Hoch^\mb(\cat(A,B,K))\to\Hoch^\mb(A)$ and
$\Hoch^\mb(\cat(A,B,K))\to\Hoch^\mb(B)$ are {\it quasi-isomorphisms
of complexes}? (They are always maps of dg Lie algebras, it is
clear). The answer is given as follows: it is enough if the
following conditions are satisfied:

Consider the left action of $B$ on $K$. It is a map of right
$A$-modules, and we get a map $L_B^0\colon B\to\Hom_{mod-A}(K,K)$.
We can also derive this map to a map $L_B\colon B\to
\RHom_{mod-A}(K,K)$. Analogously, we define from the right
$A$-action on $K$ the map $R_A\colon A^{\opp}\to\Rhom_{B-mod}(K,K)$.
\begin{definition}
Let $A$ and $B$ be two dg associative algebras, and let $K$ be dg
$B$-$A$-bimodule. We say that the triple $(A,B,K)$ is a Keller's
admissible triple if the maps
\begin{equation}\label{final2.2.1}
\begin{aligned}
\ &L_B\colon B\to \RHom_{mod-A}(K,K)\\
&R_A\colon A^{\opp}\to \RHom_{B-mod}(K,K)
\end{aligned}
\end{equation}
are quasi-isomorphisms of algebras.
\end{definition}

There are known two examples when the Keller's condition is
satisfied:
\begin{itemize}
\item[(1)] $A$ is any dg associative algebra, and there is a map
$\varphi\colon B\to A$ which is {\it a quasi-isomorphism} of
algebras. We set $K=A$ with the tautological structure of right
$A$-module on it, and with the left $B$-module structure given by
the map $\varphi$;
\item[(2)] $A$ is a quadratic Koszul algebra, $B=A^!$ is the Koszul
dual algebra, and $K$ is the Koszul complex of $A$ considered as a
$B$-$A$-bimodule.
\end{itemize}
The both statements are proven in [Kel1]. The theory developed in
Section 1 makes the generalization of (2) for Koszul algebras over
discrete valuation rings straightforward.

The following theorem was found and proven in [Kel1]:
\begin{theorem*}
Let $(A,B,K)$ be a Keller's admissible triple. Then the natural
projections $p_A\colon\Hoch^\mb(\cat(A,B,K))\to\Hoch^\mb(A)$ and
$p_B\colon\Hoch^\mb(\cat(A,B,K))\to\Hoch^\mb(B)$ are
quasi-isomorphisms of dg Lie algebras.
\begin{proof}
Let $t\colon L^\mb\to M^\mb$ be a map of complexes. Recall that its
cone $\Cone(t)$ is defined as $\Cone(t)=L^\mb[1]\oplus M^\mb$ with
the differential given by matrix
$$
d=\begin{pmatrix}d_{L[1]}& 0\\
t[1]& d_M
\end{pmatrix}
$$
To prove that the map $t\colon L^\mb\to M^\mb$ is a
quasi-isomorphism, it is equivalently than to prove that the cone
$\Cone(t)$ is acyclic in all degrees.

Let us consider the cone $\Cone(p_A)$ where $p_A\colon
\Hoch^\mb(\cat)\to\Hoch^\mb(A)$ is the natural projection. Let us
prove that if the first condition of (\ref{final2.2.1}) is
satisfied, the cone $\Cone(p_A)$ is acyclic.

We can regard $\Cone(p_A)$ as a bicomplex where the vertical
differentials are the Hochschild differentials and the horizontal
differential is $p_A[1]$. This bicomplex has two columns, therefore
its spectral sequences converge. Compute firstly the differential
$p_A[1]$. Then the term $E_1$ is the sum of $\Hoch^\mb(B)\oplus
\Hoch^\mb(B,K,A)$, as a graded vector space. There are 3 components
of the differential in $E_1$: the Hochschild differentials in
$\Hoch^\mb(B)$ and in $\Hoch^\mb(B,K,A)$, and exactly the same
differential $d_\Hoch^{BK}\colon\Hoch^\mb(B)\to\Hoch^\mb(B,K,A)[1]$,
as in the Hochschild complex of the category
$\Hoch^\mb(\cat(A,B,K))$.

Compute firstly the cohomology of $\Hoch^\mb(B,K,A)$ with the only
Hochschild differential. One can write:
\begin{equation}\label{final2.2.2}
\Hoch^\mb(B,K,A)=\Hom(T(B),\Hom(K\otimes T(A),K))
\end{equation}
with some differentials, where we denote by $T(V)$ the free
associative algebra generated by $V$. More precisely, the term
$\Hom_{\mathbb{C}}(K\otimes T(A),K)$ is equal to the complex
$\Hom_{mod-A}(\mathrm{Bar}^\mb_{mod-A}(K),K)$ of maps from the
bar-resolution of $K$ in the category of right $A$-modules to $K$.
This is equal to $\Rhom_{mod-A}(K,K)$, which is quasi-isomorphic to
$B$ by the first Keller's condition. But this is not all what we
need--we also need to know that the left $B$-module structures on
$B$ and on $\Rhom_{mod-A}(K,K)$ are the same. This is exactly
guaranteed by the Keller's condition, which says that the
quasi-isomorphism $B\to\Rhom_{mod-A}(K,K)$ {\it is induced by the
left action of $B$ on $K$}.

Now we have two complexes, which are exactly the same, and are
$\Hoch^\mb(B)$, but there is also the component $d_\Hoch^{BK}$ from
one to another. In other words, so far our complex is the cone of
the identity map from $\Hoch^\mb(B)$ to itself, and this cone is
clearly acyclic.

We have proved that if the first Keller's condition is satisfied,
the natural projection $p_A\colon \Hoch^\mb(\cat)\to\Hoch^\mb(A)$ is
a quasi-isomorphism. If the second Keller's condition is satisfied,
we conclude, analogously, that the projection $p_B\colon
\Hoch^\mb(\cat)\to\Hoch^\mb(B)$ is a quasi-isomorphism.
\end{proof}
\end{theorem*}
B. Keller used this theorem in [Kel1] to show that in the two cases
listed above when the Keller's conditions are satisfied, the
Hochschild cohomological complexes of $A$ and $B$ are
quasi-isomorphic as dg Lie algebras. In particular, this is true
when $A$ and $B$ are Koszul dual algebras, the case of the most
interest for us.

\begin{remark}
If $A$ and $B$ are Koszul dual algebras, but $K$ is replaced by
$\mathbb{C}$, the only cohomology of the Koszul complex, we still
have the quasi-isomorphisms
$B\to\Rhom_{mod-A}(\mathbb{C},\mathbb{C})$ and
$A^{\opp}\to\Rhom_{B-mod}(\mathbb{C},\mathbb{C})$, {\it but these
maps are not induced by the left (correspondingly, right) actions of
$B$ (correspondingly, $A$) on $\mathbb{C}$}. These actions define
some stupid maps which are not quasi-isomorphisms. This example
shows that the Keller's dg category in this case may be not
quasi-equivalent (and it is really the case) to its homology dg
category.
\end{remark}

\subsubsection{The Keller's condition in the (bi)graded case}
As we already mentioned in Remark in Section 0.3, when the algebras $S(V^*)$ and $\Lambda(V)$ are considered just as associative algebras, they are not Koszul dual. Namely, $\Ext_{\Lambda(V)}(k,k)=S[[V^*]]$, the formal power series instead of polynomials. 
To avoid this problem, we should work in the category of algebras with {\it inner} $\mathbb{Z}$-grading and with {\it cohomological} $\mathbb{Z}$-grading. Finally, $\Lambda^k(V)$ should have the inner grading $k$ and the cohomological grading $k$, while $S^k(V^*)$ has the inner grading $k$ and the cohomological grading 0. 
Then we should switch to the category of bigraded modules, and compute $\Ext$ algebras in this category. In this definition such $\Ext$ algebras will be automatically bigraded. This completely agrees with the theory of Koszul dually discussed in Section 1.

The only problem is that the Keller's Theorem 4.2 was proven above for the category of graded algebras when only the cohomological grading is taken into the account. To make this Theorem valid for the bigraded case, we should modify the definition of the Hochschild cohomological complex of a bigraded algebra and of a bigraded dg category. 

We give the following definition. Let $A$ be a bigraded algebra (one grading is inner and another one is cohomological). We define the graded Hochschild complex $\Hoch^\mb_{\mathrm{gr}}(A)$ as the direct sum of its bigraded components. The same definition will be done for bigraded dg categories. 

In general, it is {\bf not} true that the graded Hochschild complex is quasi-isomorphic to the usual one (which is graded only with respect to the cohomological grading but not with respect to the inner). They are quasi-isomorphic for $A=S(V^*)$, but for $A=\Lambda(V)$ this quasi-isomorphism fails. Indeed, the cohomology of the usual Hochschild complex for $A=\Lambda(V)$ is the {\it formal} polyvector fields on $V$ while the cohomology of the graded Hochschild complex is in this case the {\it polynomial} polyvector fields on $V$. Therefore, if we need to work in the bigraded category, we should reprove Theorem 4.2 in this case.

Let $A$ and $B$ be two associative bigraded algebras, and let $K$ be dg bigraded $B-A$-bimodule. We say that the triple $(A,B,K)$ satisfies the graded Keller's condition if the natural maps 
\begin{equation}\label{ref100}
\begin{aligned}
\ &L_B\colon B\to\mathrm{RHom}_{grmod-A}(K,K)\\
&R_A\colon A^{\opp}\to \RHom_{B-grmod}(K,K)
\end{aligned}
\end{equation}
are {\it bigraded} quasi-isomorphisms. Here the derived functors are taken in the categories of graded modules.

\begin{theorem*}
Let $(A,B,K)$ satisfies the graded Keller's condition. Then the natural projections $p_A\colon\Hoch^\mb_\gr(\cat(A,B,K))\to\Hoch^\mb_\gr(A)$
and $p_B\colon \Hoch^\mb_\gr(\cat(A,B,K))\to\Hoch^\mb_\gr(B)$ are quasi-isomorphisms.
\end{theorem*}
 
The proof is completely analogous to the usual case, and we leave the details to the reader.\qed

In the sequel we will omit the subscript $\gr$ with the notation of the Hochschild complex, always assuming the theory developing here.
\subsection{The maps $p_A$ and $p_B$ are maps of $B_\infty$
algebras} Let $A,B$ be two associative algebras, and let $K$ be any
$B-A$-bimodule, not necessarily satisfying the Keller's condition
from Section 4.2. Then we have two projections $p_A\colon
\Hoch^\mb(\cat(A,B,K))\to \Hoch^\mb(A)$ and $p_B\colon
\Hoch^\mb(\cat(A,B,K))\to\Hoch^\mb(B)$. We know from Section 2 that
the Hochschild complex $\Hoch^\mb(A)$ of any associative algebra has
the natural structure of $B_\infty$ algebra by means of  the
Getzler-Jones' braces (see Figure 2). The same is true for
$\Hoch^\mb(\mathcal{C})$ where $\mathcal{C}$ is a dg category, which
is established by the same braces' construction.

The following simple Lemma, due to Bernhard Keller [Kel1], is very
important for our paper:

\begin{lemma}
Let $A,B$ be two associative algebras, and let $K$ be a
$B-A$-bimodule. Then the natural projections $p_A\colon
\Hoch^\mb(\cat(A,B,K))\to \Hoch^\mb(A)$ and $p_B\colon
\Hoch^\mb(\cat(A,B,K))\to\Hoch^\mb(B)$ are maps of $B_\infty$
algebras.
\begin{proof}
It is clear because the projections $p_A$ and $p_B$ are compatible
with the braces, and with the cup-products. That is, they are
compatible with the maps $m_i$ and $m_{ij}$ of the $B_\infty$
structure, see Section 2.5.
\end{proof}
\end{lemma}

\subsection{We formulate a new version of the Main Theorem}
Let $A=S(V^*)$, $B=\Lambda(V)$, and $K=K^\mb(S(V^*))$. We know from
Section 1.1 the the algebra $S(V^*)$ is Koszul, and its Koszul dual
$A^!=\Lambda(V)$. Thus, we can apply Theorem 4.2 to the triple
$(S(V^*),\Lambda(V),K^\mb(S(V^*)))$. We have constructed a
$B_\infty$ algebra $\Hoch^\mb(\cat(A,B,K))$ for $A,B,K$ as above,
and the diagram
\begin{equation}\label{eqq1}
\xymatrix{& \Hoch^\mb(\cat(A,B,K))\ar[dl]_{p_A}\ar[dr]^{p_B}\\
\Hoch^\mb(A)&&\Hoch^\mb(B)}
\end{equation}
where the two right maps are maps of $B_\infty$ algebras. Let now
$t\colon G_\infty\to B_\infty$ be a point of the Tamarkin's
manifold, see Section 3.3. Then the diagram (\ref{eqq1}) is a
diagram of maps of $G_\infty$ algebras, depending on $t\in
\mathfrak{T}$.

Let now $\U=\mathfrak{X}^0(t)$ be the universal $G_\infty$ morphism
$\mathcal{G}_V\colon T_\poly(V)\to \Hoch^\mb(S(V^*))$ defined for
all finite-dimensional (graded) vector spaces $V$, see Section 3.4.
It depends on the point $t\in \mathfrak{T}$ and is defined up to a
homotopy. Denote by $\mathcal{G}^S(t)$ and $\mathcal{G}^\Lambda(t)$
the specializations of this universal $G_\infty$ morphism for the
vector spaces $V$ and $V^*[1]$, correspondingly. Identify
$T_\poly(V)$ with $T_\poly(V^*[1])$ as in Section 0.1 of the
Introduction. Then we have the following diagram of $G_\infty$ maps:

\begin{equation}\label{eqq2}
\xymatrix{&\Hoch^\mb(A)\\
T_\poly(V)\ar[ur]^{\mathcal{G}^S(t)}\ar[dr]_{\mathcal{G}^\Lambda(t)}&&\Hoch^\mb(\cat(A,B,K))\ar[ul]_{p_A}\ar[dl]^{p_B}\\
&\Hoch^\mb(B)}
\end{equation}
depending on $t\in \mathfrak{T}$.

Here and in Sections 6 we prove the following statement:
\begin{theorem*}
For any fixed $t\in \mathfrak{T}$, the diagram (\ref{eqq2}) is
homotopically commutative, that is, it is commutative in the
Quillen's homotopical category.
\end{theorem*}

Now restrict ourselves with the $L_\infty$ component of the
$G_\infty$ maps. Then clearly the diagram remains to be
homotopically commutative. We have the following

\begin{corollary}{\bf(A new version of the Main Theorem)}
For any $t\in \mathfrak{T}$, the diagram (\ref{eqq2}) defines a
homotopically commutative diagram of $L_\infty$ maps.
\end{corollary}

We explain in Section 7 in detail why this Corollary implies the
Main Theorem in our previous version, for Koszul duality in
deformation quantization.

Now let us begin to prove the Theorem above.

{\it Proof of Theorem (beginning)}: The proof is based on the
following Key-Lemma:
\begin{klemma}
For any $t\in \mathfrak{T}$, the diagramm (\ref{eqq2}) defines a
commutative diagram of isomorphisms maps on cohomology.
\end{klemma}
We prove this Lemma in Section 6, and it will take some work.

Now let us explain how the Theorem follows from the Key-Lemma.

The diagram (\ref{eqq2}) is a diagram of $G_\infty$
quasi-isomorphisms (the two left arrows clearly are
quasi-isomorphisms, and the two right ones are by the Keller's
Theorem proven in Section 4.2). We can uniquely up to a homotopy
invert a $G_\infty$ quasi-isomorphism. Then the map
\begin{equation}\label{4new1}
\mathcal{G}(t)=(\mathcal{G}^\Lambda(t))^{-1}\circ p_B\circ
p_A^{-1}\circ\mathcal{G}^S(t)
\end{equation}
is uniquely defined, up to a homotopy, $G_\infty$ quasi-automorphism
of $T_\poly(V)$. Now, by the Key-Lemma, its first Taylor component
is the identity map. Then we can take the logarithm
\begin{equation}\label{4new2}
\mathcal{D}=\log(\mathcal{G})
\end{equation}
which is a $G_\infty$ derivation of $T_\poly(V)$. We are in the situation of the following lemma:

\begin{lemma}
Let $\mathcal{G}$ be an $\Aff(V)$-equivariant $G_\infty$ automorphism of the Gerstenhaber algebra $T_\poly(V)$ with the standard Gerstenhaber structure, whose first component is the identity map. Then the $G_\infty$ automorphism $\mathcal{G}$ is the identity. 
\begin{proof}
As above, we can take $\mathcal{D}=\log\mathcal{G}$, then $\mathcal{D}$ is an $\Aff(V)$-invariant $G_\infty$ derivation of $T_\poly(V)$.
By Theorem 2.4, this $G_\infty$ derivation is 0. Therefore, $\mathcal{G}=\exp\mathcal{D}$ is the identity.
\end{proof}
\end{lemma}

The Theorem is now proven mod out the Key-Lemma which we prove in Section
6.

\section{The homotopical category of dg algebras over a Koszul
operad} Here we give, following [Sh3], a construction of the
homotopy category, appropriate for our needs in the next Sections of
this paper. Our emphasis here is how the homotopy relation reflects
in the gauge equivalence condition for deformation quantization. We
restrict ourselves with the case of the operad of Lie algebras
because this is the only case we will use. The constructions for
general Koszul operad are analogous.

Here we use the construction of Quillen homotopical category given
in [Sh3]. In a sense, it is "the right cylinder homotopy relation".
Recall here the definition.

\subsection{The homotopy relation from [Sh3]} Let $\g_1,\g_2$ be
two dg Lie algebras. Then there is a dg Lie algebra
$\Bbbk(\g_1,\g_2)$ which is pro-nilpotent and such that the
solutions of the Maurer-Cartan equation in $\Bbbk(\g_1,\g_2)^1$ are
exactly the $L_\infty$ morphisms from $\g_1$ to $\g_2$. Then the
zero degree component $\Bbbk(\g_1,\g_2)^0$ acts on the Maurer-Cartan
solutions, as usual in deformation theory (the dg Lie algebra
$\Bbbk(\g_1,\g_2)$ is pro-nilpotent), and this action gives a
homotopy relation.

The dg Lie algebra $\Bbbk(\g_1,\g_2)$ is constructed as follows. As
a dg vector space, it is
\begin{equation}\label{hf2.40}
\Bbbk(\g_1,\g_2)=\Hom(C_+(\g_1,\mathbb{C}),\g_2)
\end{equation}
Here $C(\g_1,\mathbb{C})$ is the chain complex of the dg Lie algebra
$\g_1$, it is naturally a counital dg coalgebra, and
$C_+(\g_1,\mathbb{C})$ is the kernel of the counit map.

Define now a Lie bracket on $\Bbbk(\g_1,\g_2)$. Let
$\theta_1,\theta_2\in\Bbbk(\g_1,\g_2)$ be two elements. Their
bracket $[\theta_1,\theta_2]$ is defined (up to a sign) as
\begin{equation}\label{eq2.1}
C_+(\g_1,\mathbb{C})\xrightarrow{\Delta}C_+(\g_1,\mathbb{C})^{\otimes
2}\xrightarrow{\theta_1\otimes\theta_2}\g_2\otimes\g_2\xrightarrow{[,]}g_2
\end{equation}
where $\Delta$ is the coproduct in $C_+(\g_1,\mathbb{C})$ and $[,]$
is the Lie bracket in $\g_2$. It follows from the cocommutativity of
$\Delta$ that in this way we get a Lie algebra.

An element $F$ of degree 1 in $\Bbbk(\g_1,\g_2)$ is a collection of
maps
\begin{equation}\label{eq2.2}
\begin{aligned}
\ & F_1\colon \g_1\to\g_2\\
&F_2\colon \Lambda^2(\g_1)\to\g_2[-1]\\
&F_3\colon \Lambda^3(\g_1)\to\g_2[-2]\\
&\dots
\end{aligned}
\end{equation}
and the Maurer-Cartan equation $d_{\Bbbk}F+\frac12[F,F]_{\Bbbk}=0$
is the same that the collection $\{F_i\}$ are the Taylor components
of an $L_\infty$ map which we denote also by $F$. Note that the
differential in $\Bbbk(\g_1,\g_2)$ comes from 3 differentials: the
both inner differentials in $\g_1$ and $\g_2$, and from the chain
differential in $C_+(\g_1,\mathbb{C})$.

Now the solutions of the Maurer-Cartan equation form a quadric in
$\g^1$, and for any pro-nilpotent dg Lie algebra $\g$, the component
$\g^0$ acts on (the pro-nilpotent completion of) this quadric by
vector fields. Namely, each $X\in\g^0$ defines a vector field
\begin{equation}\label{eq2.3}
\frac{dF}{dt}=-dX+[X,F]
\end{equation}
It can be directly checked that this vector field indeed preserves
the quadric.

In our case, this vector field can be exponentiated to an action on
the pro-nilpotent completion on $\Bbbk$. This action gives our
homotopy relation on $L_\infty$ morphisms.

\subsection{Application to deformation quantization} Let
$\g_1,\g_2$ be two dg Lie algebras, and let $\F^1,\F^2\colon \g_1\to
\g_2$ be two homotopic in the sense of Section 2.4.1 $L_\infty$
morphisms.

Let $\alpha$ be a solution of the Maurer-Cartan equation in $\g_1$.
Any $L_\infty$ morphism $\F\colon\g_1\to \g_2$ gives a solution
$\F_*\alpha$ of the Maurer-Cartan equation in $\g_2$, by formula

\begin{equation}\label{hf2.50}
\F_*\alpha=\F_1(\alpha)+\frac12\F_2(\alpha\wedge\alpha)+\frac16\F_3(\alpha\wedge\alpha\wedge\alpha)+\dots
\end{equation}
(suppose that this infinite sum makes sense).

Then in our situation we have two solutions $\F^1_*\alpha$ and
$\F^2_*\alpha$ of the Maurer-Cartan equation in $\g^2$.

\begin{lemma}
Suppose that all infinite sums (exponents) we need make sense in our
situation. Suppose two $L_\infty$ morphisms $\F_1,\F_2\colon
\g_1\to\g_2$ are homotopic in the sense of Section 2.4.1, and
suppose that $\alpha$ is a solution of the Maurer-Cartan equation in
$\g_1$. Then the two solutions $\F^1_*\alpha$ and $\F^2_*\alpha$ of
the Maurer-Cartan equation in $\g_2$ are gauge equivalent.
\begin{proof}
Let $X\in \Bbbk(\g_1,\g_2)^0$ be the generator of the homotopy
between $\F^1$ and $\F^2$. Define
\begin{equation}\label{hf2.51}
X_*\alpha=X(\alpha)+\frac12X(\alpha\wedge\alpha)+\frac16X(\alpha\wedge\alpha\wedge\alpha)+\dots
\end{equation}
Then $X_*\alpha\in (\g_2)^0$. Consider the vector field on
$(\g_2)^1$:
\begin{equation}\label{hf2.52}
\frac{dg}{dt}=-d(X_*\alpha)+[X_*\alpha,g]
\end{equation}
Then the exponent of this vector field maps $\F^1*\alpha$ to
$\F^2_*\alpha$.
\end{proof}
\end{lemma}

\section{The main computation}
Here we prove the Key-Lemma 4.4 which is only remains to conclude
the proof of Theorem 4.4.
\subsection{}
We are going to construct "the Hochschild-Kostant-Rosenberg map"
$\varphi_{HKR}^\cat\colon T_\poly(V)\to\Hoch^\mb(\cat(A,B,K))$ where
$A=S(V^*)$, $B=\Lambda(V)$, and $K$ is the Koszul complex of
$S(V^*)$. At the final step of the compuation, we normalize the
Koszul differential by $\dim V$, as follows:
\begin{equation}\label{6.1.1}
d_\Koszul^{\mathrm{norm}}=\frac1{\dim V}\sum_{a=1}^{\dim
V}x_a\frac{\partial}{\partial \xi_a}
\end{equation}
However, in the computation below we suppose that the Koszul complex
is not normalized. The normalized Koszul complex defines the
equivalent Keller's category, so it is irrelevant.

Our Hochschild-Kostant-Rosenberg map $\varphi_{HKR}^\cat$ will make
the following diagram commutative (up to a sign) on the cohomology:

\begin{equation}\label{6.1.2}
\xymatrix{&\Hoch^\mb(A)\\
T_\poly(V)\ar[ur]^{\mathcal{G}^S(t)}\ar[dr]_{\mathcal{G}^\Lambda(t)}\ar@{.>}[rr]^{\varphi_{HKR}^\cat}&&\Hoch^\mb(\cat(A,B,K))\ar[ul]_{p_A}\ar[dl]^{p_B}\\
&\Hoch^\mb(B)}
\end{equation}
We did not specify the sign, but it does not make any problem.

In the computation below we use the graphical representation of the
cochains in $\Hom(\Lambda(V)^{\otimes m}\otimes K\otimes S(V^*),K)$.
The reader familiar with the Kontsevich's paper [K97] will
immediately understand our graphical representation. (But for other
readers, we define our cochains by the explicit formulas, see
(\ref{eqhkr1})-(\ref{eqhkr3}) below).

In our graphical cochains, we consider a circle with two fixed
points, 0 and $\infty$. The arguments from $\Lambda(V)$ are placed
on the left half of the circle, and the arguments from $S(V^*)$ are
placed on the right half. Any arrow is the operator
$\sum_{a=1}^{\dim V}\frac{\partial}{\partial \xi_a}\cdot
\frac{\partial}{\partial x_a}$. In our convention, which coincides
with the one in [K97], the start-point of any arrow "differentiates"
the odd arguments, while the end-point differentiates the even
arguments. We have one point inside the disc bounded by the circle,
where we place our polyvector field $\gamma$. We use the notation
$\gamma=\gamma^S\otimes \gamma^\Lambda$ (where $\gamma^S$ and
$\gamma^\Lambda$ are the even and the odd coordinates of $\gamma$)
and suppose that $\gamma$ is homogeneous in the both $x_i$'s and
$\xi_j$'s coordinates.

After this general remarks, let us start.

\subsection{Some graph-complex}
The problem of a construction of quasi-isomorphism
$\varphi_{HKR}^\cat\colon T_\poly(V)\to\Hoch^\mb(\cat)$ is rather
non-trivial. Indeed, the usual Hochschild-Kostant-Rosenberg cochains
$\varphi_{HKR}^S(\gamma)\in Hoch(S(V^*))$ and
$\varphi_{HKR}^\Lambda(\gamma)\in\Hoch^\mb(\Lambda(V))$ are not
cocycles when considered as cochains in $\Hoch^\mb(\cat)$. Indeed,
their boundaries have components which belong in $\Hom(K^\mb\otimes
S(V^*)^{\otimes m_1},K^\mb)$ and in $\Hom(\Lambda(V)^{\otimes
m_1}\otimes K^\mb, K^\mb)$, correspondingly. Our map
$\varphi_{HKR}^\cat$ contain as summand the both cochains
$\varphi_{HKR}^S$ and $\varphi_{HKR}^\Lambda$, and many other
summands. These other summands are the cochains associated with the
graphs $F_{0,m_2}^0$ and $F_{m_1,0}^\infty$ shown in Figure below.

\sevafigc{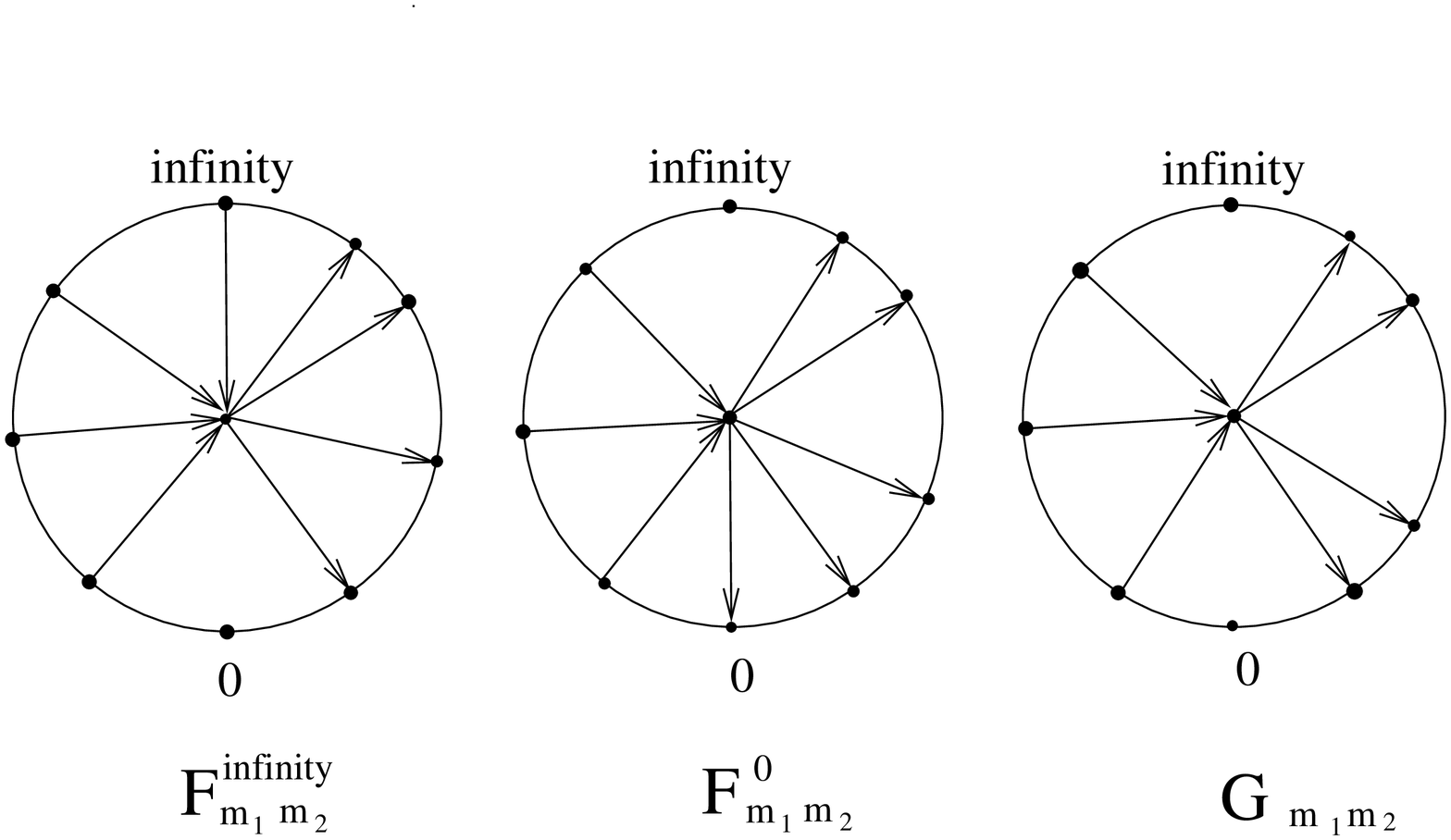}{120mm}{0}{The cochains
$F_{m_1,m_2}^\infty$, $F_{m_1,m_2}^0$, and $G_{m_1,m_2}$ for $m_1=3,
m_2=4$}

It is instructive to formulate the following Proposition in a bit
more generality than we really need, for all graphs $F^0_{m_1,m_2}$
and $F^\infty_{m_1,m_2}$. Denote the corresponding maps
$\Phi_\Gamma$ in $\Hom(\Lambda(V)^{\otimes m_1}\otimes K^\mb\otimes
S(V^*)^{\otimes m_2}, K^\mb)$ by $F_{m_1,m_2}^\infty(\gamma)$,
$F_{m_1,m_2}^0(\gamma)$, and $G_{m_1,m_2}(\gamma)$, where $\gamma\in
T_\poly(V)$. Suppose that $\gamma$ is homogeneous in both $x_i$'s
and $\xi$'s. As maps $T_\poly(V)\to \Hoch^\mb(\cat)$ the maps
$F_{m_1,m_2}^0$ and $F_{m_1,m_2}^\infty$ have degree 0, and the map
$G_{m_1,m_2}$ has degree 1.

We have the following explicit formulas for these maps:
\begin{equation}\label{eqhkr1}
\begin{aligned}
\
&G_{m_1,m_2}(\gamma)(\lambda)=\frac1{n!}\frac1{m!}\sum_{i_1,..,i_{m_1}=1}^{\dim
V}\sum_{j_1,...,j_{m_2}=1}^{\dim V} \\&\pm k\left(\lambda\wedge
\partial{\xi_{i_1}}(\lambda_1)\wedge\dots\wedge
\partial{\xi_{i_{m_1}}}(\lambda_{m_1})\wedge
(\partial\xi_{j_1}\circ\dots\circ\partial\xi_{j_{m_2}}(\gamma^\Lambda))\right)\times\\
&\times(\partial x_{i_1}\circ\dots\circ\partial
x_{i_{m_1}})(\gamma^S)\cdot
\partial x_{j_1}(f_1)\dots\partial x_{j_{m_2}}(f_{m_2})
\end{aligned}
\end{equation}

\begin{equation}\label{eqhkr2}
\begin{aligned}
\
&F_{m_1,m_2}^0(\gamma)(\lambda)=\frac1{m_1!}{(m_2+1)!}\sum_{i_1,..,i_{m_1}=1}^{\dim
V}\sum_{a,j_1,...,j_{m_2}=1}^{\dim V}\\ & \pm {\mathbf{\partial
x_a}}(k)\left(\lambda\wedge
\partial{\xi_{i_1}}(\lambda_1)\wedge\dots\wedge
\partial{\xi_{i_{m_1}}}(\lambda_{m_1})\wedge
({\mathbf{\partial \xi_a}}\circ\partial\xi_{j_1}\circ\dots\circ\partial\xi_{j_{m_2}}(\gamma^\Lambda))\right)\times\\
&\times(\partial x_{i_1}\circ\dots\circ\partial
x_{i_{m_1}})(\gamma^S)\cdot
\partial x_{j_1}(f_1)\dots\partial x_{j_{m_2}}(f_{m_2})
\end{aligned}
\end{equation}

\begin{equation}\label{eqhkr3}
\begin{aligned}
\ &F_{m_1,m_2}^\infty(\gamma)(\lambda)=
\frac1{(m_1+1)!}\frac1{m_2!}\sum_{b,i_1,..,i_{m_1}=1}^{\dim
V}\sum_{j_1,...,j_{m_2}=1}^{\dim V}\\ &\pm
k\left(\mathbf{\partial\xi_b}\left(\lambda\wedge
\partial{\xi_{i_1}}(\lambda_1)\wedge\dots\wedge
\partial{\xi_{i_{m_1}}}(\lambda_{m_1})\wedge
(\partial\xi_{j_1}\circ\dots\circ\partial\xi_{j_{m_2}}(\gamma^\Lambda))\right)\right)\times\\
&\times(\mathbf{\partial x_b}\circ\partial
x_{i_1}\circ\dots\circ\partial x_{i_{m_1}})(\gamma^S)\cdot
\partial x_{j_1}(f_1)\dots\partial x_{j_{m_2}}(f_{m_2})
\end{aligned}
\end{equation}
Here, as usual, we denote by $\{x_i\}$ some basis in $V^*$, and by
$\{\xi_i\}$ the dual basis in $V[-1]$.

Let $\gamma$ be a polynomial polyvector field in $T_\poly(V)$,
homogeneous in both $x$'s and $\xi$'s. Denote $\deg_S\gamma$ and
$\deg_\Lambda\gamma$ the corresponding homogeneity degrees. (We have
the Lie degree $\deg\gamma=\deg_\Lambda\gamma-1$).

Denote $d_\Hoch$ and $d_\Koszul$ the Hochschild and Koszul
components of the differential acting on $\Hoch^\mb(\Lambda(V),
K^\mb, S(V^*))\subset \Hoch^\mb(\cat)$.

\begin{proposition}
Suppose $\sharp In_\Gamma(v)\le \deg_S\gamma$ and $\sharp Star(v)\le
\deg_\Lambda\gamma$ for each separate graph $\Gamma$ in the claims
below, where $v$ is the only vertex of the first type. Suppose that
$F_{m,n}^0(\gamma)$ etc. means the sum over all orderings of the
sets $Star(v)$ and $In(v)$ (see (10) and (11) in the definition of
an admissible graph), that is, over all admissible graphs which are
the same geometrically. (The sum should be taken with the
appropriate signs depending naturally on the orderings). Then we
have:
\begin{itemize}
\item[(i)] $d_\Hoch F_{m,n}^0(\gamma)=\pm G_{m, n+1}(\gamma)$,
\item[(ii)] $d_\Koszul F_{m,n}^0(\gamma)=\pm \dim V\cdot(\deg_\Lambda(\gamma)-n)\cdot G_{m,n}(\gamma)$,
\item[(iii)] $d_\Hoch F_{m,n}^\infty(\gamma)=\pm G_{m+1,n}(\gamma)$,
\item[(iv)] $d_\Koszul F_{m,n}^\infty(\gamma)=\pm \dim V\cdot (\deg_S(\gamma)-m)\cdot  G_{m,n}(\gamma)$.
\end{itemize}
\begin{proof}
The proof of Proposition is just a straightforward computation. For
convenience of the reader, we present it here in all details.

We give the proofs of (i) and (ii); the proofs of the second two
statements are analogous.

Prove (i).

It would be instructive for the reader to recall before the proof
the proof that the classical Hochschild-Kostant-Rosenberg
$\varphi_{HKR}^S(\gamma)$ is a Hochschild cocycle in
$\Hoch^\mb(S(V^*))$ for any $\gamma\in T_\poly(V)$. It goes as
follows: we associate with a $k$-polyvector field $\gamma$ the
cochain $\varphi_{HKR}^S(\gamma)\in \Hom(S(V^*)^{\otimes k},
S(V^*))$ defined as
\begin{equation}\label{eqhkr4}
\varphi_{HKR}^S(\gamma)(f_1\otimes\dots\otimes
f_k)=\sum_{i_1,...,i_k=1}^{\dim V}\pm
\gamma(dx_{i_1}\wedge\dots\wedge dx_{i_k})\partial
x_{i_1}(f_1)\dots\partial x_{i_k}(f_k)
\end{equation}
The only nonzero terms may appear when all $i_1,\dots,i_k$ are
different,and the sign $\pm$ is the sign of the permutation
$(1,2,\dots,k)\mapsto(i_1,i_2,\dots, i_k)$. The proof that
$\varphi_{HKR}^S(\gamma)$ is a Hochschild cocycle just uses the
Leibniz formula $\partial x_a(f_i f_{i+1})=\partial
x_a(f_i)f_{i+1}+f_i\partial x_a(f_{i+1})$ and the Hochschild
coboundary formula
\begin{equation}\label{eqhkr5}
\begin{aligned}
\ &d_\Hoch(\Psi)(f_1\otimes\dots\otimes f_{k+1})=\\
&f_1\Psi(f_2\otimes f_2\otimes\dots\otimes f_{k+1})-\\
&-\Psi((f_1f_2)\otimes f_3\otimes\dots)+\Psi(f_1\otimes
(f_2f_3)\otimes\dots)\mp\dots\\
&+(-1)^{k+1}\Psi(f_1\otimes\dots\otimes f_k)f_{k+1}
\end{aligned}
\end{equation}
We see that the all terms will be mutually canceled. Now let us see
when this kind of phenomenon may be destroyed in the coboundary of
$F_{m,n}^0(\gamma)$. It is clear that any problem place is the
marked point 0 at the boundary of the circle. Consider the sum of
two "problematic" summands. This is
\begin{equation}\label{eqhkr6}
\pm
F_{m,n}^0(\gamma)(\lambda_1\otimes\dots\otimes\lambda_{m}\otimes(\lambda_{m+1}(k))\otimes
f_n\otimes\dots\otimes f_1)\mp
F_{m,n}^0(\gamma)(\lambda_1\otimes\dots\otimes
\lambda_m\otimes((k)f_{n+1})\otimes f_{n}\otimes\dots\otimes f_1)
\end{equation}
Here we use the notation $\lambda(k)$ and $(k)f$ for the left action
of $\Lambda(V)$ and for the right action of $S(V^*)$,
correspondingly. These two summands give from (\ref{eqhkr2})
\begin{equation}\label{eqhkr7}
\pm \partial x_a(\lambda_{m+1}(k))=\pm\lambda_{m+1}(\partial x_a k)
\end{equation}
which clearly is canceled with (a part of) the previous summand,
\begin{equation}\label{eqhkr8}
\mp
F_{m,n}^0(\gamma)(\lambda_1\otimes\dots\otimes(\lambda_{m}\lambda_{m+1})\otimes
k\otimes f_n\otimes\dots\otimes f_1)
\end{equation}
So the first summand in (\ref{eqhkr6}) does not contribute to the
answer. Contrary, the second summand gives the term
\begin{equation}\label{eqhkr9}
\partial x_a(k\cdot f_{n+1})=\partial x_a(k)\cdot f_{n+1}+k\cdot \partial
x_a(f_{n+1})
\end{equation}
The first summand in (\ref{eqhkr9}) is canceled with the one of two
summands in $F_{m,n}^0(\gamma)(\lambda_1\otimes\dots\otimes
\lambda_m\otimes k\otimes(f_{n+1}\cdot f_{n})\otimes\dots\otimes
f_1)$. The second summand in (\ref{eqhkr9}) is not canceled with an
other summand, and it gives the only term which contributes to the
answer. This term clearly gives $G_{m, n+1}(\gamma)$.

Prove (ii).

We need to compute
\begin{equation}\label{eqhkr10}
\begin{aligned}
\ &\sum_{i_1,..,i_{m_1}=1}^{\dim V}\sum_{a,j_1,...,j_{m_2}=1}^{\dim
V} \pm {\mathbf{\large d_\Koszul}}\\
&\{{\mathbf{\partial x_a}}{\mathbf{\large (k)}}\left(\lambda\wedge
\partial{\xi_{i_1}}(\lambda_1)\wedge\dots\wedge
\partial{\xi_{i_{m_1}}}(\lambda_{m_1})\wedge
({\mathbf{\partial \xi_a}}\circ\partial\xi_{j_1}\circ\dots\circ\partial\xi_{j_{m_2}}(\gamma^\Lambda))\right)\times\\
&\times(\partial x_{i_1}\circ\dots\circ\partial
x_{i_{m_1}})(\gamma^S)\cdot
\partial x_{j_1}(f_1)\dots\partial x_{j_{m_2}}(f_{m_2})\}\\
&\mp \sum_{i_1,..,i_{m_1}=1}^{\dim
V}\sum_{a,j_1,...,j_{m_2}=1}^{\dim V} \\
&\pm {\mathbf{\partial x_a}}{\mathbf{\large{(d_\Koszul
k)}}}\left(\lambda\wedge
\partial{\xi_{i_1}}(\lambda_1)\wedge\dots\wedge
\partial{\xi_{i_{m_1}}}(\lambda_{m_1})\wedge
({\mathbf{\partial \xi_a}}\circ\partial\xi_{j_1}\circ\dots\circ\partial\xi_{j_{m_2}}(\gamma^\Lambda))\right)\times\\
&\times(\partial x_{i_1}\circ\dots\circ\partial
x_{i_{m_1}})(\gamma^S)\cdot
\partial x_{j_1}(f_1)\dots\partial x_{j_{m_2}}(f_{m_2})
\end{aligned}
\end{equation}
We have:
\begin{equation}\label{eqhkr11}
d_\Koszul k=\sum_{p=1}^{\dim V}x_p\partial \xi_p
\end{equation}
Then (\ref{eqhkr10}) is equal to
\begin{equation}\label{eqhkr12}
\begin{aligned}
\ &\sum_{i_1,..,i_{m_1}=1}^{\dim V}\sum_{a,j_1,...,j_{m_2}=1}^{\dim
V}\\
& \pm {\mathbf{\large x_p\partial\xi_p}}\{{\mathbf{\partial
x_a}}{\mathbf{\large (k)}}\left(\lambda\wedge
\partial{\xi_{i_1}}(\lambda_1)\wedge\dots\wedge
\partial{\xi_{i_{m_1}}}(\lambda_{m_1})\wedge
({\mathbf{\partial \xi_a}}\circ\partial\xi_{j_1}\circ\dots\circ\partial\xi_{j_{m_2}}(\gamma^\Lambda))\right)\times\\
&\times(\partial x_{i_1}\circ\dots\circ\partial
x_{i_{m_1}})(\gamma^S)\cdot
\partial x_{j_1}(f_1)\dots\partial x_{j_{m_2}}(f_{m_2})\}\\
&\mp \sum_{i_1,..,i_{m_1}=1}^{\dim
V}\sum_{a,j_1,...,j_{m_2}=1}^{\dim V} \\
&\pm {\mathbf{\partial x_a}}{\mathbf{\large{(x_p\partial\xi_p
k)}}}\left(\lambda\wedge
\partial{\xi_{i_1}}(\lambda_1)\wedge\dots\wedge
\partial{\xi_{i_{m_1}}}(\lambda_{m_1})\wedge
({\mathbf{\partial \xi_a}}\circ\partial\xi_{j_1}\circ\dots\circ\partial\xi_{j_{m_2}}(\gamma^\Lambda))\right)\times\\
&\times(\partial x_{i_1}\circ\dots\circ\partial
x_{i_{m_1}})(\gamma^S)\cdot
\partial x_{j_1}(f_1)\dots\partial x_{j_{m_2}}(f_{m_2})
\end{aligned}
\end{equation}
where the summation over $p$ is assumed. Clearly (\ref{eqhkr12}) is
$A+B$ where
\begin{equation}\label{eqhkr13}
\begin{aligned}
\ &{\mathbf A}=\\
&\sum_{i_1,..,i_{m_1}=1}^{\dim V}\sum_{a,j_1,...,j_{m_2}=1}^{\dim
V}\\
&\pm {\mathbf{\large x_p}}{\mathbf\Large
[}{{\partial\xi_p}}\{{\mathbf{\partial x_a}}{\mathbf{\large
(k)}}\left(\lambda\wedge
\partial{\xi_{i_1}}(\lambda_1)\wedge\dots\wedge
\partial{\xi_{i_{m_1}}}(\lambda_{m_1})\wedge
({\mathbf{\partial \xi_a}}\circ\partial\xi_{j_1}\circ\dots\circ\partial\xi_{j_{m_2}}(\gamma^\Lambda))\right)\times\\
&\times(\partial x_{i_1}\circ\dots\circ\partial
x_{i_{m_1}})(\gamma^S)\cdot
\partial x_{j_1}(f_1)\dots\partial x_{j_{m_2}}(f_{m_2})\}\\
&\mp \sum_{i_1,..,i_{m_1}=1}^{\dim
V}\sum_{a,j_1,...,j_{m_2}=1}^{\dim V}\\
&\pm {\mathbf{\large{(\partial\xi_p \partial x_a
k)}}}\left(\lambda\wedge
\partial{\xi_{i_1}}(\lambda_1)\wedge\dots\wedge
\partial{\xi_{i_{m_1}}}(\lambda_{m_1})\wedge
({\mathbf{\partial \xi_a}}\circ\partial\xi_{j_1}\circ\dots\circ\partial\xi_{j_{m_2}}(\gamma^\Lambda))\right)\times\\
&\times(\partial x_{i_1}\circ\dots\circ\partial
x_{i_{m_1}})(\gamma^S)\cdot
\partial x_{j_1}(f_1)\dots\partial
x_{j_{m_2}}(f_{m_2}){\mathbf\Large ]}
\end{aligned}
\end{equation}
and
\begin{equation}\label{eqhkr14}
\begin{aligned}
\ &{\mathbf{B}}=\\
&\mp {\mathbf{\dim V}}\sum_{i_1,..,i_{m_1}=1}^{\dim
V}\sum_{a,j_1,...,j_{m_2}=1}^{\dim V}\\
&\pm {\mathbf{\large{(\partial\xi_a k)}}}\left(\lambda\wedge
\partial{\xi_{i_1}}(\lambda_1)\wedge\dots\wedge
\partial{\xi_{i_{m_1}}}(\lambda_{m_1})\wedge
({\mathbf{\partial \xi_a}}\circ\partial\xi_{j_1}\circ\dots\circ\partial\xi_{j_{m_2}}(\gamma^\Lambda))\right)\times\\
&\times(\partial x_{i_1}\circ\dots\circ\partial
x_{i_{m_1}})(\gamma^S)\cdot
\partial x_{j_1}(f_1)\dots\partial
x_{j_{m_2}}(f_{m_2})
\end{aligned}
\end{equation}
In the last equation the symbol $\delta_{ap}$ appears when we take
the commutator $[\partial x_a, x_p]$ in the second summand of
(\ref{eqhkr12}), which gives the factor $\dim V$ and summation only
over $a$ in $B$. We continue for $A$ and $B$ separately.

Let us start with $B$. We have:
\begin{equation}\label{eqhkr15}
\begin{aligned}
\ &B=\\
&\mp {\dim V}\sum_{i_1,..,i_{m_1}=1}^{\dim
V}\sum_{a,j_1,...,j_{m_2}=1}^{\dim V} \\
&\pm {\mathbf{k}}\left(\lambda\wedge
\partial{\xi_{i_1}}(\lambda_1)\wedge\dots\wedge
\partial{\xi_{i_{m_1}}}(\lambda_{m_1})\wedge
({\mathbf{\xi_a\partial \xi_a}}\circ\partial\xi_{j_1}\circ\dots\circ\partial\xi_{j_{m_2}}(\gamma^\Lambda))\right)\times\\
&\times(\partial x_{i_1}\circ\dots\circ\partial
x_{i_{m_1}})(\gamma^S)\cdot
\partial x_{j_1}(f_1)\dots\partial
x_{j_{m_2}}(f_{m_2})=\\
&\mp {\dim V\cdot
(\deg_\Lambda(\gamma)-m_2)}\times\\
&\times\sum_{i_1,..,i_{m_1}=1}^{\dim
V}\sum_{a,j_1,...,j_{m_2}=1}^{\dim V}\\
&\pm {\mathbf{k}}\left(\lambda\wedge
\partial{\xi_{i_1}}(\lambda_1)\wedge\dots\wedge
\partial{\xi_{i_{m_1}}}(\lambda_{m_1})\wedge
\partial\xi_{j_1}\circ\dots\circ\partial\xi_{j_{m_2}}(\gamma^\Lambda))\right)\times\\
&\times(\partial x_{i_1}\circ\dots\circ\partial
x_{i_{m_1}})(\gamma^S)\cdot
\partial x_{j_1}(f_1)\dots\partial
x_{j_{m_2}}(f_{m_2})=\\
&=\pm \dim V\cdot (\deg_\Lambda(\gamma)-m_2)\cdot G_{m,n}(\gamma)
\end{aligned}
\end{equation}

Now turn back to the computation of $A$. Clearly (up to the sign,
but the signs always work for us) that $A=0$. Indeed, schematically
the formula (\ref{eqhkr13}) for $A$ looks like
$\partial\xi_p(k(\lambda\wedge T))-(\partial\xi_p(k))(\lambda\wedge
T)$ for some $T\in\Lambda(V)$. If we define
$k^\prime(\lambda)=k(\lambda\wedge T)$ we need to compute
\begin{equation}\label{eqhkr16}
(\partial\xi_p(k^\prime))(\lambda)-(\partial\xi_p(k))(\lambda\wedge
T) \end{equation} But
$(\partial\xi_p(k^\prime))(\lambda)=k^\prime(\xi_p\wedge\lambda)=k(\xi_p\wedge\lambda\wedge
T)$. Now we see that the two summands in (\ref{eqhkr16}) are equal.

We have proved the statements (i) and (ii) of the Proposition. The
proofs of (iii) and (iv) are analogous.
\end{proof}
\end{proposition}

\subsection{Construction of the Hochschild-Kostant-Rosenberg map
$\varphi^{\cat}_{HKR}$} Now we have everything we need to construct
the map $\varphi^{\cat}_{HKR}\colon T_\poly(V)\to\Hoch^\mb(\cat)$.
Of course, it would be better to specify the signs in the
Proposition above; however, we will see that the construction below
does not depend seriously on these signs.

Suppose that $\deg_S\gamma=m$, $\deg_\Lambda\gamma=n$.  Start with
$\wtilde{\varphi}^S_{HKR}(\gamma)\in\Hoch^\mb(S(V^*))$, which is by
definition the Hochschild-Kopstant-Rosenberg cochain without
division by the $n!$. It total differential in $\Hoch^\mb(\cat)$ is
$d_\tot\wtilde{\varphi}_{HKR}^S(\gamma)=(\pm)G_{0,n}(\gamma)$. From
now on, we will suppose that the all signs in Proposition above are
$"+"$, if some of them are $"-"$, the formula will be the same up to
some signs. So suppose that
$d_\tot\wtilde{\varphi}_{HKR}^S(\gamma)=G_{0,n}(\gamma)$ with sign
$+$. We know from statement (i) of the Proposition that $d_\Hoch
F_{0,n-1}^0(\gamma)=G_{à,n}(\gamma)$, the same cochain. Therefore,
$d_\Hoch(\wtilde{\varphi}_{HKR}^S(\gamma)-F_{0,n-1}^0(\gamma))=0$.
But then $\wtilde{\varphi}_{HKR}^S(\gamma)-F_{0,n-1}^0(\gamma)$ has
a non-trivial Koszul differential which can be found by Proposition
(ii). We have:
$d_\Koszul(\wtilde{\varphi}_{HKR}^S(\gamma)-F_{0,n-1}^0(\gamma))=d_\Koszul(F_{0,n-1}^0(\gamma))=
\dim V\cdot G_{0,n-1}(\gamma)$. Now we want to kill this coboundary
by the Hochschild differential. We have: $d_\Hoch(\dim V\cdot
F^0_{0,n-2}(\gamma))=\dim V\cdot G_{0,n-1}(\gamma)$. Continuing in
this way, we find that (we omit $\gamma$ at each term):
\begin{equation}\label{eqhkrfinal1}
\begin{aligned}
\ &d_\tot(\wtilde{\varphi}_{HKR}^S-F^0_{0,n-1}+\dim V\cdot
F^0_{0,n-2}-\dots+\dots+(-1)^n(n-1)!\dim^{n-1}V\cdot F^0_{0,0})\\
&=(-1)^n n!\dim^nVG_{0,0}
\end{aligned}
\end{equation}
But we can start also with $\wtilde{\varphi}_{HKR}^\Lambda(\gamma)$,
and finally get also $G_{0,0}$ with some multiplicity. More
precisely, we have:
\begin{equation}\label{eqhkrfinal2}
\begin{aligned}
\ &d_\tot(\wtilde{\varphi}_{HKR}^\Lambda-F^\infty_{m-1,0}+\dim
V\cdot F^\infty_{m-2,0}-\dots+(-1)^m(m-1)!\dim^{m-1}V\cdot
F^\infty_{0,0})\\
&=(-1)^m m!\dim^m V G_{0,0}
\end{aligned}
\end{equation}
We finally set:
\begin{equation}\label{eqhkrfinal}
\begin{aligned}
\varphi_{HKR}^\cat=&(-1)^n\frac1{n!\dim^n
V}\left(\wtilde{\varphi}_{HKR}^S+\sum_{i=1}^n(-1)^i
(i-1)!\dim^{i-1}V
F^0_{0,n-i}\right)-\\
&(-1)^m\frac1{m!\dim^m
V}\left(\wtilde{\varphi}_{HKR}^\Lambda+\sum_{j=1}^m(-1)^j(j-1)!\dim^{j-1}V
F^\infty_{m-j,0}\right)
\end{aligned}
\end{equation}
It is a cocycle in the Hochschild cohomological complex
$\Hoch^\mb(\cat)$:
\begin{equation}\label{eqhkr}
d_\tot\varphi_{HKR}^\cat(\gamma)=0
\end{equation}
for any $\gamma\in T_\poly(V)$.

We can prove the following
\begin{theorem*}
The map $\varphi_{HKR}^\cat\colon T_\poly(V)\to\Hoch^\mb(\cat)$ is a
quasi-isomorphism of complexes. When we use the normalized Koszul
differential instead of the usual one (so, it has the same effect as
to set $\dim V=1$ in the formula above), the map
$\varphi_{HKR}^\cat$ makes the diagram (\ref{6.1.2}) commutative (up
to a non-essential sign) on the level of cohomology.
\begin{proof}
The second statement is clear. The first one (that
$\varphi_{HKR}^\cat$ is a quasi-isomorphism of complexes) follows
from the second one and from Theorem 4.2 which says that the maps
$p_A$ and $p_B$ are quasi-isomorphisms in our case.
\end{proof}
\end{theorem*}

Key-Lemma 4.4 is proven.\qed

Theorem 4.4 is proven.\qed

\section{Proof of the Main Theorem}
First of all, we formulate the Main Theorem exactly in the form we
will prove it here.
\subsection{The final formulation of the Main Theorem}
\begin{theorem*}{\bf (Main Theorem, final form)}
Suppose $t\colon G_\infty\to B_\infty$ is a quasi-isomorphism of
operads, and let $\U_V=\mathfrak{X}(t)_V\colon
T_\poly(V)\to\Hoch^\mb(S(V^*))$ be the corresponding $L_\infty$ map,
defined uniquely up to homotopy (see Section 3). Let $\alpha$ be a
quadratic Poisson bivector on $V$, and let $\mathcal{D}(\alpha)$ be
the corresponding quadratic Poisson bivector on $V^*[1]$. Denote by
$S(V^*)_\hbar$ and $\Lambda(V)_\hbar$ the corresponding deformation
quantizations of $S(V^*)\otimes\mathbb{C}[[\hbar]]$ and
$\Lambda(V)\otimes\mathbb{C}[[\hbar]]$ given by
\begin{equation}\label{7.1.1}
f\star g=f\cdot g+\hbar\cdot \U_1(\alpha)(f\otimes
g)+\frac12\hbar^2\cdot \U_2(\alpha\wedge\alpha)(f\otimes g)+\dots
\end{equation}
Then the algebras $S(V^*)_\hbar$ and $\Lambda(V)_\hbar$ are graded
(where $\deg\hbar=0$, $deg x_i=1$ for all $i$) and quadratic. Also,
they are Koszul as algebras over the discrete valuation ring
$\mathbb{C}[[\hbar]]$, see Section 1.  Moreover, they are Koszul
dual to each other.
\end{theorem*}

We prove the Theorem throughout this Section.

\subsection{An elementary Lemma}
We start with the following simple statement:
\begin{lemma}
\begin{itemize}
\item[(1)] Suppose $K_\hbar$ is a free $\mathbb{C}[[\hbar]]$-module,
which is also a left (or right) $\mathbb{C}[[\hbar]]$-linear module
over an algebra $A_\hbar$ which is supposed to be also free as
$\mathbb{C}[[\hbar]]$-module. Then if the specialization
$K_{\hbar=0}$ is a free module over the specialization
$A_{\hbar=0}$, $K_\hbar$ is a free left (right) $A_\hbar$-module;
\item[(2)] suppose $K^\mb_\hbar$ is a complex of free
$\mathbb{C}[[\hbar]]$-modules ($\deg\hbar=0$) with
$\mathbb{C}[[\hbar]]$-linear differential. Suppose that the $i$-th
cohomology (for some $i$) of the specialization $K^\mb_{\hbar=0}$ is
zero. Then the $i$-th cohomology of $K_\hbar^\mb$ is also zero.
\end{itemize}
\begin{proof}
The both statements are standard; let us recall the proofs for
convenience of the reader.

(1): Suppose the contrary, then for some $k_i(\hbar)\in K_\hbar$ and
some $a_i(\hbar)\in A_\hbar$ one has $\sum_ia_i(\hbar)\cdot
k_i(\hbar)=0$. Let $N$ be the minimal power of $\hbar$ in the
equation. Then we can divide the equation over $\hbar^N$ and the
equation still holds, because the both $A_\hbar$ and $K_\hbar$ are
free $\mathbb{C}[[\hbar]]$-modules. Then we reduce over $\hbar$ and
get a nontrivial linear equation for the $A_{\hbar=0}$-module
$K_{\hbar=0}$ which contradicts to the assumption.

(2): Let $k_i(\hbar)$ be an $i$-cicycle in $K_\hbar^\mb$, we should
prove that it is a coboundary. Suppose $\hbar^N$ is the minimal
power of $\hbar$ in $k_i(\hbar)$, then we divide over $\hbar^N$. We
get again a cocycle, because the differential is
$\mathbb{C}[[\hbar]]$-linear and $K^\mb_\hbar$ is a free
$\mathbb{C}[[\hbar]]$-module. Denote this new cocycle again by
$k_i(\hbar)$. Then its zero degree in $\hbar$ term is a cocycle in
the reduced complex $K_{\hbar=0}$ and we can kill it by some
coboundary. Then substract and divide over minimal power of $\hbar$,
ans so on.
\end{proof}
\end{lemma}

\subsection{The algebras $S(V^*)_\hbar$ and $\Lambda(V)_\hbar$ are
Koszul} We start to prove the Theorem. Prove firstly that the
algebras $S(V^*)_\hbar$ and $\Lambda(V)_\hbar$ are graded quadratic
and Koszul. The first statement is proven analogously to the
speculation in Section 0.2. The difference that here in a universal
deformation quantization we may have more general graphs than in the
Kontsevich's quantization, namely non-connected graphs and graphs
with simple loops. But it does not change the proof.

Let us prove that these algebras are Koszul. Consider the case of
$S(V^*)_\hbar$, the proof for $\Lambda(V)_\hbar$ is analogous.

By Lemma 1.2.5, it is necessary to prove that the Koszul complex
$K_\hbar^\mb=\bigl(S(V^*)_\hbar\otimes_{\mathbb{C}[[\hbar]]}\Hom_{\mathbb{C}[[\hbar]]}(S(V^*)^!,\mathbb{C}[[\hbar]]),d_\Koszul\bigr)$
is acyclic in all degrees except degree 0. The complex $K_\hbar^\mb$
is clearly a complex of free $\mathbb{C}[[\hbar]]$-modules with a
$\mathbb{C}[[\hbar]]$-linear differential. We are in situation of
Lemma 7.2(2), because the specialization at $\hbar=0$ gives clearly
the Koszul complex for the usual algebra $S(V^*)$ which is known to
be acyclic. We are done.

\subsection{We continue to prove the Main Theorem}
Now we prove the only non-trivial part of the Theorem, that the
algebras $S(V^*)_\hbar$ and $\Lambda(V)_\hbar$ are Koszul dual.

Consider the diagram (\ref{eqq2}). It is a diagram of $G_\infty$
quasi-isomorphisms which is known to be homotopically commutative,
see Theorem 4.4. Then we can construct a $G_\infty$
quasi-isomorphism $\F\colon T_\poly(V)\to\Hoch^\mb(\cat(A,B,K))$
dividing the diagram into two commutative triangles. Restrict $\F$
to its $L_\infty$ part. Then we get an $L_\infty$ quasi-isomorphism
$\F\colon T_\poly(V)\colon \Hoch^\mb(\cat(A,B,K))$. Here
$A=S(V^*)\otimes \mathbb{C}[[\hbar]]$, $B=\Lambda(V)\otimes
\mathbb{C}[[\hbar]]$, etc.

Then this $L_\infty$ map $\F$ attaches to the Maurer-Cartan solution
$\alpha\in T_\poly(V)$ (our quadratic Poisson bivector field) a
solution of the Maurer-Cartan equation in $\Hoch^\mb(\cat(A,B,K))$,
by formula
\begin{equation}\label{7.4.1}
\F_*(\alpha)=\hbar\F_1(\alpha)+\frac12\hbar^2\F_2(\alpha\wedge\alpha)+\dots
\end{equation}

What a solution of the Maurer-Cartan equation in
$\Hoch^\mb(\cat(A,B,K))$ means im more direct terms?

It consists from the following data:

\begin{itemize}
\item[(i)] A deformation quantization $A_\hbar$ of the algebra
$A=S(V^*)\otimes \mathbb{C}[[\hbar]]$;
\item[(ii)] a deformation quantization $B_\hbar$ of the algebra
$B=\Lambda(V)\otimes\mathbb{C}[[\hbar]]$;
\item[(iii)] a deformed differential on the Koszul complex
$K^\mb(S(V^*))\otimes \mathbb{C}[[\hbar]]$, we denote the deformed
complex by $K_\hbar^\mb$;
\item[(iv)] a structure of a $B_\hbar$-$A_\hbar$-bimodule on
$K_\hbar^\mb$.
\end{itemize}

The crucial point is the following Lemma:

\begin{lemma}
The algebra $A_\hbar$ is gauge equivalent (and therefore isomorphic)
to the algebra $S(V^*)_\hbar$ from Section 7.1, and the algebra
$B_\hbar$ is gauge equivalent to $\Lambda(V)_\hbar$.
\begin{proof}
It follows from the commutativity of the diagram (\ref{eqq2}), and
from Lemma 5.2.
\end{proof}
\end{lemma}

\subsection{We finish to prove the Main Theorem}
From Lemma 7.4, it is enough to prove that the quadratic graded
algebras $A_\hbar$ and $B_\hbar$ are Koszul dual to each other. For
this (because the both algebras are Koszul) it is enough to prove
that $B_\hbar=A_\hbar^!$. Let us prove it.

The complex $K_\hbar$ is a complex of $B_\hbar$-$A_\hbar$ modules.
As complex of $A_\hbar$-modules, it is free by Lemma 7.2(1). By
Lemma 7.2(2), it is a free $A_\hbar$-resolution of the module
$\mathbb{C}[[\hbar]]$. Therefore, we can use $K_\hbar$ for the
computation of the Koszul dual algebra:
\begin{equation}\label{7final.1}
(A_\hbar)^!=\mathrm{RHom}_{Mod-A_\hbar}(K_\hbar, K_\hbar)
\end{equation}
On the other hand, from the bimodule structure (see (iv) in the list
in Section 7.4), we have an {\it algebra} homomorphism
\begin{equation}\label{7final.2}
B_\hbar\rightarrow \mathrm{RHom}_{Mod-A_\hbar}(K_\hbar, K_\hbar)
\end{equation}
We only need to prove that it is an isomorphism. It again follows
from the facts that the both sides are free
$\mathbb{C}[[\hbar]]$-modules (for the l.h.s. it is clear, for the
r.h.s. it follows from (\ref{7final.1})), and that the
specialization of (\ref{7final.2}) at $\hbar=0$ is an isomorphism.

Theorem 7.1 is proven.\qed

Faculty of Science, Technology and Communication, Campus
Limpertsberg, University of Luxembourg,
162A avenue de la Faiencerie, L-1511 LUXEMBOURG\\
{\it e-mail}: {\tt borya$\_$port@yahoo.com}

\end{document}